\documentclass[english]{article}
\usepackage[T1]{fontenc}
\usepackage[utf8]{inputenc}
\usepackage{morefloats}
\usepackage[letterpaper]{geometry}
\usepackage{float}
\usepackage{amsmath}

\usepackage{graphicx}
\usepackage{color}
\usepackage{epstopdf}
\usepackage{adjustbox}
\usepackage{diagbox}
\usepackage{multirow}
\usepackage{array}
\usepackage{mathrsfs}
\usepackage{amsmath, amssymb, amsthm}
\usepackage{subfigure}
\usepackage{tikz}

\graphicspath{{./Figures/}}

\newtheorem{definition}{Definition}
\newtheorem{theorem}{Theorem}
\newtheorem{assumption}{Assumption}
\newtheorem{lemma}{Lemma}

\usepackage{amsfonts} 
\usepackage{color}
\usepackage{adjustbox}
\usepackage{diagbox}
\definecolor{black}{rgb}{0,0,0}

\definecolor{red}{rgb}{1,0,0}

\definecolor{blue}{rgb}{0,0,1}

\usepackage{multirow}

\makeatother

\usepackage{babel}
\usepackage{authblk}
\title{}

\title{\textbf{A local-global generalized multiscale finite element method for
		highly heterogeneous stochastic groundwater flow problems}}
\author{Yiran Wang\thanks{Department of Mathematics, The Chinese University of Hong Kong, Shatin, Hong Kong SAR.}, \;
	Eric Chung\thanks{Department of Mathematics, The Chinese University of Hong Kong, Shatin, Hong Kong SAR.} \;
	and \; Shubin Fu\thanks{Department of Mathematics, University of Wisconsin-Madison, WI, USA. Corresponding author (shubinfu89@gmail.com)}
}

\begin{document}
	\maketitle
	\begin{abstract}
		In this paper, we propose a local-global multiscale method for highly heterogeneous stochastic groundwater flow problems under the framework of reduced basis method and 
		the generalized multiscale finite element method (GMsFEM).
		Due to incomplete characterization of the medium properties of the groundwater flow problems, random variables are used to parameterize the uncertainty. As a result, solving the problem repeatedly is required to obtain statistical quantities. 
		Besides, the medium properties are usually highly heterogeneous, which will result in a large linear system that needs to be solved. 
		Therefore, it is intrinsically inevitable to seek a computational-efficient model reduction method to overcome the difficulty. We will explore the combination of the 
		reduced basis method and the GMsFEM. In particular, we will use residual-driven basis functions, which are key ingredients in GMsFEM. This local-global multiscale method is more efficient than applying the GMsFEM or reduced basis method individually. We first construct parameter-independent multiscale basis functions that include both local and global
		information of the permeability  fields, and then use these basis functions to construct several global snapshots and global basis functions for fast online computation with different parameter inputs.
		We provide rigorous analysis of the proposed method and extensive numerical examples to demonstrate the accuracy and efficiency of the local-global multiscale method.
	\end{abstract}
	\section{Introduction}
	
	In groundwater flow problems, it is usually difficult to know the exact permeability at all points in the practical simulation domain due to large scope of the domain and inevitable natural small scale randomness. One way to deal with the  uncertainty  of the permeability
	field is to write it as a stochastic function, then the groundwater flow problems become stochastic partial differential equations (SPDEs). 
	To obtain statistics of interested quantities,  SPDEs have to be solved many times with different input parameters. 
	For the inverse problems such as estimating the permeability \cite{chen2006data,shi2012multiscale},
	one also needs to solve the underlying forward model many times to obtain convergence.
	On the other hand, the description of permeability field can be detailed at multiple scales  from pore scales to geological scales thanks to the development of geostatistical modeling techniques and characterization methods.
	As a result, direct computation in the finest scale is difficult due to limited computational resources. 
	Therefore it is necessary to develop efficient multiscale methods and model reduction
	methods to cope with above mentioned difficulties in solving heterogeneous stochastic  underground flow problems.
	
	Many multiscale methods were proposed to deal with highly heterogeneous problems,
	for example, the multiscale finite element method \cite{hou1997multiscale,chen2003mixed}, 
	variational multiscale method \cite{hughes1998variational}, multiscale mortar methods \cite{arbogast2007multiscale,wheeler2012multiscale} and multiscale finite volume methods \cite{cortinovis2014iterative,lunati2004multi}.
	Among all, the multiscale finite element method (MsFEM) and its extension 
	the Generalized Multiscale Finite element method (GMsFEM) have achieved
	huge success for various types of heterogeneous problems \cite{chung2010reduced,chung2016adaptive,efendiev2009multiscale,chung2017online,vasilyeva2019constrained,vasilyeva2019multiscale}.
	Underlying idea of the MsFEM is to  construct multiscale basis functions with appropriate designed local problems and boundary conditions, these basis functions 
	incorporate small scale medium information thus can yield accurate low dimensional coarse-grid solutions. 
	Besides, once constructed, the multiscale basis functions could be applied to cases with different source terms. However, 
	one basis function  in each coarse neighborhood is not sufficient \cite{efendiev2013generalized,chung2014adaptive} if there are long channels and non-separable scales in the permeability field. Then, the GMsFEM is proposed to handle arbitrary complex media  by constructing multiple basis functions. This method is divided into two major steps, constructing the snapshot space and then the multiscale basis function space. First, a rich snapshot space is constructed by solving a static problem with 
	multiple boundary conditions or sources, this snapshot space include all local information of the target solution. Furthermore, a reduced space named offline space is obtained by some well-designed spectral problems. 
	It is shown that when the number of offline basis functions exceeds some level, the approximation error will decrease slowly \cite{efendiev2013generalized}. In particular, the convergence rate is proportional to $1/\Lambda$, where $\Lambda$ is the smallest eigenvalue  abandoned during the construction of multiscale space. 
	To further improve the accuracy, residual-driven basis functions \cite{chung2015residual,yang2018residual,efendiev2016online} are proposed as enrichment upon previous offline space. 
	These residual basis have three major advantages. First of all, those residual-driven basis are more powerful in reducing error with only a few basis functions since they contain global residual information. In addition, one can flexibly adjust the number of such basis functions in each local neighborhood \cite{efendiev2016online,chung2018fast} and the enrichment times. Lastly, similar to previous offline computation, the construction of residual-driven basis can also be conducted in parallel due to the independent computation in different local neighborhoods. Consequently, it is computational-efficient to include residual-driven basis in the approximation space.
	
	However, for stochastic permeability field case, special treatment should be considered when using the GMsFEM. Previous efforts can be seen \cite{jiang2016reduced,he2019reduced,ma2018local,jiang2017model,chung2018cluster,chung2019adaptiveJCAM}. All these methods show high accuracy and efficiency by designing parameter-independent multiscale basis functions. However, those multiscale basis functions contain only local medium information which limits the accuracy of the methods. Our goal here is also to construct parameter-independent basis functions that can capture global information of the media by including the residual-driven basis functions in training stage.
	However, we will not instantly use these parameter-independent basis functions  for online 
	coarse-grid simulation with different input parameters to get statistical quantities. Instead, we select some sample parameters and use the parameter-independent multiscale   basis functions to get some fine-grid defined solutions, who are less accurate compared with fine-grid defined solutions obtained by FEM but still remains high-accuracy due to the global information contained in the multiscale basis functions.
	We call these fine-grid defined solutions as global snapshots and we then use the proper orthogonal decomposition (POD) to extract dominant modes for fast online 
	simulations. In particular, we utilize a special-designed spectral problem which is derived from convergence analysis to reduce the dimension of approximation space. Therefore, this is a local-global method since the multiscale basis functions are constructed   locally in the first stage and then be used to generate
	globally defined snapshots, we call this method as GMsFEM-POD method. The global snapshots here are generated in coarse grids without solving any fine-grid large-scale problems and thus more efficient compared with the standard strategies in POD method. 
	The online stage here is the same as the standard POD method and thus enjoys 
	high efficiency.  Rigorous  analysis is provided and extensive numerical experiments will be presented to verify the theoretical analysis and 
	show the superior performance of our method. In particular, we  consider 
	highly heterogeneous random medium characterized  by the Karhunen-Lo$\grave{\mathbf{e}}$ve (KL) expansion \cite{huang2001convergence}.
	We will study the influence of using different parameters in KL expansion and the influence of number of POD basis and local basis functions. We also compare our local-global method with standard POD method.  Numerical results show that 
	including residual-driven basis functions can improve the accuracy of proposed method compared with only using equal number of local offline basis functions. Besides, the accuracy of our method is comparable to the standard POD method. 
	
	The paper is organized as follows. In section \ref{sec:GMsFEM}, we briefly review the GMsFEM including  construction of offline and residual-driven basis. GMsFEM-POD method is illustrated in section 3. In section 4, we analyze  our GMsFEM-POD method and numerical results are presented in section 5.
	\section{Multiscale model reduction using the GMsFEM}\label{sec:GMsFEM}
	In this section, we present the GMsFEM for the groundwater flow problem. First, we present some preliminaries in Section 2.1. The constructions of snapshot space and offline space are presented in Section 2.2. The local online enrichment process is introduced in Section 2.3. 
	\subsection{Preliminaries}
	Let $\Omega$ be a bounded domain in $\mathrm{R^d}(d=2,3)$ and $T>0$ be a fixed time. Let $\Omega_r$ be  a sample space. The transient flow system in heterogeneous random porous media can be described as   
	\begin{eqnarray}
		\begin{aligned}
			S\dfrac{\partial u}{\partial t} - \text{div}( \kappa(\mathbf{x},\omega) \nabla u )  = & f(t,\mathbf{x}),& \quad t \in 
			[0,T], \quad \mathbf{x}\in \Omega, \quad \omega \in \Omega_r\\
			u(0,\mathbf{x})  = &g(\mathbf{x}),&  \quad \mathbf{x}\in  \Omega,\\
			u(t,\mathbf{x})  = &0,& \quad   t\in [0,T], \quad \mathbf{x} \in \partial \Omega,
		\end{aligned}
		\label{eqn:model}
	\end{eqnarray}
	where $u(t,\mathbf{x})$ is the pressure head, $f(t,\mathbf{x})$
	is the source term, $S$ is the specific storage which is assumed to be $1$ in following for simplicity. $\kappa(\mathbf{x},\omega)$ is a possible highly heterogeneous random permeability field. We use the usual notations such as $H^1(\Omega)$ and $H_0^1(\Omega)$ for Sobolev space. Moreover, we  define  the following inner products and norms:
	\begin{eqnarray}
		\langle u,v\rangle=\int_{\Omega} uv \text{  and  } \mathcal{A}(u,v)=\int_{\Omega} \kappa(\mathbf{x},\omega) \nabla u\cdot \nabla v .\\
		\|u\|_{L^2}^2=\int_{\Omega} u^2 ,\quad
		|u|_{H_0^1}^2=\int_{\Omega} |\nabla u|^2,\text{  and  }
		\|u\|_{a}^2=\mathcal{A}(u,u),\label{notation_norm}
	\end{eqnarray}
	where $\|\cdot\|_a$ is called energy norm.
We also define the space $H_{\kappa}(\Omega)$ as 
	\begin{align*}
		H_{\kappa}(\Omega)=\left\{ u|\int_{\Omega}\kappa \nabla u\cdot \nabla u \leq \infty\right\}.
	\end{align*}
	
	 Moreover, we recall the following Poincar$\acute{e}$ inequality,
	\begin{eqnarray}
		\|u\|_{L^2}^2\leq Q \|u\|_a^2  \quad\forall u\in H_0^1(\Omega),\label{elliptic_L2}
	\end{eqnarray}
	for some constant $Q>0$.
	
	We will first the consider the deterministic case in this section and we use $\kappa:=\kappa(\mathbf{x},\omega)$. The variational formulation for the problem (\ref{eqn:model}) is finding $u\in H^1_0(\Omega)$ such that
	\begin{eqnarray}
		\begin{aligned}
			\left\langle \dfrac{\partial u}{\partial t}, v\right\rangle+\mathcal{A}(u,v)&=\left\langle f, v\right\rangle,  \quad &t \in [0,T], \quad \mathbf{x}\in \Omega, \quad\forall v \in H_0^1(\Omega),\\
			u(0,\mathbf{x})  &= g(\mathbf{x}), \quad &\mathbf{x}\in \Omega,\\
			u(t,\mathbf{x})  &= 0, \quad &  t\in [0,T], \quad \mathbf{x} \in \partial \Omega .
		\end{aligned}
		\label{eqn:model_va}
	\end{eqnarray}
	In order to discretize system (\ref{eqn:model_va}) in time,
	we apply the implicit Euler scheme with time step $\Delta t>0$ and obtain the following
	discretization for each time $t_i=i\Delta t,i=1,2,\cdots, N_t$, where $T=N_t\Delta t$.
	\begin{eqnarray*}
		\left\langle \cfrac{u^n-u^{n-1}}{\Delta t}, v\right\rangle+\mathcal{A}(u^n,v)=\left\langle f^n, v\right\rangle,
	\end{eqnarray*}
	where $u^i=u(t_i,\mathbf{x})$ and $f^i=f(t_i,\mathbf{x})$, for $i=0,1,\cdots, N_t$.
	
\begin{figure}[ht]
		\centering
		\subfigure
		{ \includegraphics[width=0.4\textwidth]{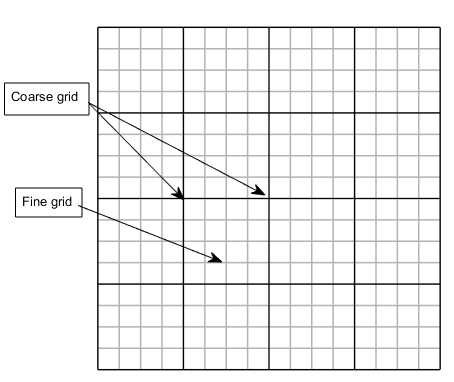}}
		\subfigure
		{ \includegraphics[width=0.33\textwidth]{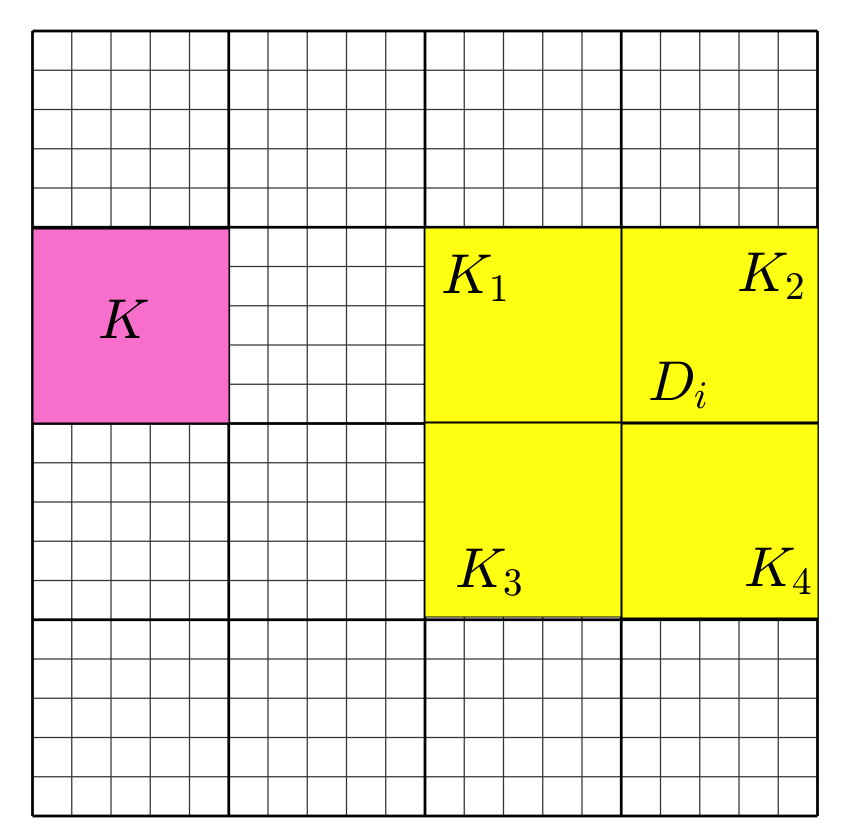}}
		
		\caption{Left: an illustration of fine and coarse grids. Right: an illustration of a coarse neighborhood and a coarse element.}\label{figure:coarse}
	\end{figure}

	We then introduce some notations about the two-scale mesh we will use.
	Let $\mathcal{T}^h$ be a finite element conforming fine partition of the domain $\Omega$, where $h>0$ is the fine grid mesh size. The coarse partition,  $\mathcal{T}^H$ of the domain $\Omega$, is formed such that each element in  $\mathcal{T}^H$ is a connected union of fine-grid elements. More precisely, $\forall K_{j} \in \mathcal{T}^H$, $ K_{j}=\bigcup_{F\in I_{j} }F$ for some $I_{j}\subset \mathcal{T}^h$, where $H>0$ is the coarse mesh size. We  consider the rectangular coarse elements and the methodology can be used with general coarse elements. We denote the interior  nodes of $\mathcal{T}^H$ by $x_i,i=1,\cdots,N_{\text{in}}$,
	where $N_\text{in}$ is the number of interior nodes. The coarse elements
	of $\mathcal{T}^H$ are denoted by $K_j,j=1,2,\cdots,N_e$, where $N_e$ is the number of coarse elements. Moreover, we define the coarse neighborhoods of the nodes $x_i$ by $D_i:=\cup\{K_j\in \mathcal{T}^H:x_i\in \overline{K_j}\}$.\\
	An illustration of the mesh  is shown in the Figure \ref{figure:coarse}, which includes fine and coarse grids, coarse elements and neighborhoods. We  suppress the time and spatial variables $t$ and $\mathbf{x}$ in functions when no ambiguity occurs.
	\subsection{The GMsFEM and the multiscale basis functions}
	Let $V_h$ be the standard finite element space spanned by piecewise linear functions $\gamma_1,\ldots,\gamma_{m}$ defined on fine mesh $\mathcal{T}^h$, where $m$ is the number of interior fine-grid nodes. Then the fine-grid solution $u_h$ is obtained in $V_h$ by solving the following problem:
	\begin{eqnarray}
		\begin{aligned}
			\frac{1}{\Delta t}\left\langle u_h^n, v\right\rangle+\mathcal{A}\left(u_h^n, v\right)&=\left\langle \frac{1}{\Delta t} u_{h}^{n-1}+f^n, v\right\rangle, \quad \forall v \in V_h,\\
			\left\langle u_{h}^{0}, v\right\rangle&=\left\langle g, v\right\rangle.
		\end{aligned}
		\label{fine sol}
	\end{eqnarray}
	 To obtain a smaller finite dimensional approximation space $V_{ms}$, the multiscale space, one needs to follow two general steps. First,  constructing a set of local snapshot basis functions in order to incorporate all  possible modes for the solutions. We note that these local snapshots are supported in local neighborhood, which are different and should be distinguished with the global snapshots  in Section \ref{sect: POD}. In the second step, we seek multiscale basis functions with a suitable spectral problem defined in the snapshot space. We take the first few dominated eigenfunctions as basis functions and use them to obtain a reduced approximation of $u$.
	
	To obtain the multiscale basis functions, we first define the snapshot space. For each coarse neighborhood $ D_{i}$, define $J_h (D_{i})$ as the set of the fine nodes of $\mathcal{T}^{h}$ lying on $\partial D_{i}$ and denote its cardinality by $L_i \in \mathbb{N}^{+}$. For each fine-grid node $x_j \in J_h( D_{i})$, we define a fine-grid function $\delta_{j}^{h}$ on $J_h( D_{i})$ as $\delta_{j}^{h}(\mathbf{x}_k)=\delta_{j,k}$. Here
	$\delta_{j,k}=1$ if $j=k$ and $\delta_{j,k}=0$ if $j\neq k$. For each $j=1,\cdots, L_i$, we define the snapshot $\psi_{j}^{(i)}$ ($j=1,\cdots,L_i$) (supported in $D_i$) as the solution to the following system:\\
	\begin{eqnarray}
		\begin{aligned}
			-\text{div}\left(\kappa \nabla \psi_{j}^{(i)}\right) &=0 \quad &\text { in } D_{i}, \\
			\psi_{j}^{(i)} &=\delta_{j}^{h} \quad &\text { on } \partial D_{i}.\label{snap_basis}
		\end{aligned}
	\end{eqnarray}
	The local snapshot space $V_{snap}^{(i)}$ corresponding to the coarse neighborhood $ D_{i}$  is defined as follows
	$V_{snap}^{(i)}:=$ \text{span}$\{\psi_{j}^{(i)}:j=1,\cdots,L_{i}\}$ and the snapshot space reads $V_{snap} :=\bigoplus_{i=1}^{N_{\text {in}}} V_{snap}^{(i)}$.
	
	In the second step, a dimension reduction is performed on each $V_{snap}^{i}$ (for $i=1,\cdots, N_{\text {in}}$) via a spectral decomposition.
	Specifically, for each $i=1,\cdots, N_{\text {in}}$, we solve the following spectral problems:
	\begin{eqnarray}
		\int_{D_{i}} \kappa \nabla \phi_{j}^{(i)} \cdot \nabla v=\lambda_{j}^{(i)} \int_{D_{i}} \hat{\kappa} \phi_{j}^{(i)} v, \quad \forall v \in V_{snap}^{(i)}, \quad j=1, \ldots, L_{i},\label{eigen}
	\end{eqnarray}
	where $\hat{\kappa} :=\kappa \sum_{i=1}^{N_{i n}} H^{2}\left|\nabla \chi_{i}\right|^{2}$ and $\{\chi_{i}\}_{i=1}^{N_{i n}}$ is a set of partition of unity that solves the following system:
	\begin{eqnarray}
		\begin{array}
			{rlrl}{-\nabla \cdot\left(\kappa \nabla \chi_{i}\right)} & {=0} & {} & {\text { in } K \subset D_{i}}, \\ {\chi_{i}} & {=p_{i}} & {} & {\text { on each } \partial K \text { with } K \subset D_{i}}, \\ {\chi_{i}} & {=0} & {} & {\text { on } \partial D_{i}}, \label{POU}
		\end{array}
	\end{eqnarray}
	where $p_i$ are some polynomial functions and we can choose linear functions for simplicity.
	Assume that the eigenvalues obtained from (\ref{eigen}) are arranged in ascending order and we may use the first $l_i$ (with $1<l_i \leq L_{i},\quad l_{i} \in
	\mathbb{N}^{+}$) eigenfunctions (related to the smallest $l_i$ eigenvalues) to
	form the local multiscale space $V_{\text{ms}}^{(i)}:=$ snap$\{\chi_{i}\phi_{j}^{(i)}:j=1,\cdots,l_{i}\}$. The mulitiscale space $V_{\text{ms}}^{(i)}$ is the direct sum of the local mulitiscale spaces, namely
	$V_{ms} :=\bigoplus_{i=1}^{N_{\text {in}}} V_{ms}^{(i)}$. 
	Moreover we define
	\begin{align}
		\Lambda=\min_{1\leq i\leq N_{in}} \lambda_{l_i+1}^{(i)}.\label{tot_eigen_val}
	\end{align} 
	Let $M$ be the total number of multiscale basis, then we write $V_{ms}=\{\zeta_1,\ldots,\zeta_M\}$.
	Once the multiscale space $V_{ms}$ is constructed, we can find the
	GMsFEM solution $u_H^n\in V_{ms}$ at $t=t_n$ by solving the following equation
	\begin{eqnarray}
		\begin{aligned}
			\frac{1}{\Delta t}\left\langle u_H^n, v\right\rangle+ \mathcal{A}\left(u_H^n, v\right)&=\left\langle \frac{1}{\Delta t} u_{H}^{n-1}+f^{n}, v\right\rangle, \quad\forall v \in V_{ms}, \\
			\langle u_{H}^{0},v\rangle&=\langle g,v\rangle, \quad\forall v \in V_{ms}.\label{ms sol}
		\end{aligned}
	\end{eqnarray}
	To derive matrix form of (\ref{fine sol}) and (\ref{ms sol}), 
	we define $c_{h}^{n}, c_{H}^{n}\in \mathbf{R}^m$ to be the coefficient vectors of $u_h^n, u_H^n\in V_h$ respectively. In particular, we have 
	\begin{eqnarray}
		u_h^n=\sum_{i=1}^{m}c_{h}^{n}[i]\gamma_i,
		\quad
		u_H^n=\sum_{i=1}^{m}c_{H}^{n}[i]\gamma_i,
	\end{eqnarray}
	where $c_{h}^{n}[i],c_{H}^{n}[i]\gamma_i$ are the $i$-th entries of $c_{h}^{n}$ and $c_{H}^{n}$ respectively. Let $\zeta_{i,c}\in \mathbf{R}^m$ be the coefficient vectors to multiscale basis $\zeta_{i}\in V_h$, for $1\leq i \leq M$.
	We further let $V_{ms,c}=[\zeta_{1,c},\ldots,\zeta_{M,c}]\in \mathrm{R}^{m\times M}$ be the matrix form for $V_{ms}$. We use $\tilde{U}_H^n\in\mathrm{R}^M$ to denote the coefficient vector of $u_H^n$ w.r.t multiscale basis $\{\zeta_1,\ldots,\zeta_M\}$. In particular, we have 
	\begin{eqnarray}
		u_H^n=\sum_{i=1}^{M}\tilde{U}_H^n[i]\zeta_i.
	\end{eqnarray}
	Moreover, we define matrices $A,M\in \mathrm{R}^{m\times m}$ as follows, where $A,M$ are known as the stiffness and mass matrices.
	\begin{eqnarray}
	   A_{i,j}=\int_{\Omega}\kappa(x)\nabla \gamma_i\cdot \nabla\gamma_j,\quad
		M_{i,j}=\int_{\Omega} \gamma_i \gamma_j.
		\label{stiff}
	\end{eqnarray}
	Let $F\in \mathrm{R}^m $ be the load vector with $F_i=\int_{\Omega}f\gamma_i$. 
	Then we have the following matrix formulation for (\ref{fine sol}).
	\begin{eqnarray*}
		\displaystyle M\dfrac{c_{h}^{n}-c_{h}^{n-1}}{\Delta t}+Ac_{h}^{n}=F.
	\end{eqnarray*}
	We can write (\ref{ms sol}) in following way:
	\begin{eqnarray*}
		\displaystyle V_{ms,c}^TMV_{ms,c}\dfrac{\tilde{U}_{H}^{n}-\tilde{U}_{H}^{n-1}}{\Delta t}+V_{ms,c}^TAV_{ms,c}\tilde{U}_{H}^{n}=V_{ms,c}^TF.
	\end{eqnarray*}
	Finally we project the multiscale solution on fine grid to get $c_{H}^{n}$ using $c_{H}^{n}=V_{ms,c}\tilde{U}_{H}^{n}$.
	\subsection{Residual-driven basis construction}\label{sect:enrichment}
	We  present the construction of residual-driven basis functions \cite{chung2015residual} in this subsection.
	One needs to construct some additional  multiscale basis functions that contain global information in each local neighborhood.
	Let $u_H^n \in V_{ms}$ be the solution obtained in (\ref{ms sol}). Given a coarse neighborhood $D_i$, we define
	$V_i:=H_0^1(D_i)\cap V_{snap}$ equipped with the norm
	$\|v\|_{V_i}^{2}:=\int_{D_i}\kappa |\nabla {v}|^2$. We also define the global residual operator $R^n: V_{snap}\rightarrow \mathbb{R}$ and local residual operator $R_i^n: V_i\rightarrow \mathbb{R}$ by
	\begin{eqnarray}
		\mathcal{R}^{n}\left(v ; u_H^n\right) :=\int_{\Omega}\left(\frac{1}{\Delta t} u_{H}^{n-1}+f^{n}\right) v-\int_{\Omega}\left(\kappa \nabla u_H^n \cdot \nabla v+\frac{1}{\Delta t} u_H^n v\right), \quad \forall v \in V_{snap}. \label{glo res}
	\end{eqnarray}
	\begin{eqnarray}
		\mathcal{R}_{i}^{n}\left(v ; u_H^n\right) :=\int_{D_{i}}\left(\frac{1}{\Delta t} u_{H}^{n-1}+f^{n}\right) v-\int_{D_{i}}\left(\kappa \nabla u_H^n \cdot \nabla v+\frac{1}{\Delta t} u_H^n v\right), \quad \forall v \in V_{i}. \label{loc res}
	\end{eqnarray}
	The operator norms of $R^n$ and $R_i^n$ are denoted by $\|R^n\|$ and $\|R_i^n\|$, which give a measure of the size of the residuals.  The residual-driven basis functions are computed during the time-marching process, contrary to
	the offline basis functions that are pre-computed.
	
	Suppose one needs to add one new residual-driven basis function $\phi$ into the space $V_i$. The analysis in \cite{chung2015residual} suggests that the required residual-driven basis $\beta\in V_i$ is the solution to the following equation
	\begin{eqnarray}
		\mathcal{A}(\beta, v)=\mathcal{R}_{i}^{n}\left(v ; u_{H}^{n,p}\right), \quad \forall v \in V_{i}.\label{online}
	\end{eqnarray}
	We refer to $p \in \mathbb{N}$ as the level of the enrichment
	and denote the solution of (\ref{ms sol})  at time $t_n$ by $u_{H}^{n,p}$. We further define $V_{ms}^{n, p}$ to be the multiscale offline space at time $t_n$ after p-th enrichment. In particular,
	let $\mathcal{I} \subset\left\{1,2, \ldots, N_{i n}\right\}$ be the index set over some non-overlapping coarse neighborhoods. For each $i\in \mathcal{I}$, we obtain a residual-driven basis function $\beta_i\in V_i$ by solving (\ref{online}) and define $\beta_{i}^{n,p}\in V_h$ to be the global extension of $\beta_i$, i.e., $\beta_{i}=\beta_{i}^{n,p}$ in $D_i$ and vanishes in $\Omega\setminus D_i$.
	We further define $V_{ms}^{n, p+1}=V_{ms}^{n, p} \oplus \operatorname{span}\left\{\beta_{i}^{n,p} : i \in \mathcal{I}\right\}$. Moreover, $V_{ms}^{n,0}:=V_{ms}$ for time level $n\in \mathbb{N}$.

	After that, solve (\ref{ms sol}) in $V_{ms}^{n, p+1}$.
	
	We summarize the two algorithms to obtain residual-driven basis functions  be  in Table \ref{online algo1} and  Table \ref{online algo2} accordingly.
	In view of the fact that enrichment is conducted in a sequence of time steps, the difference of these two algorithms exists in the dimension of multiscale space at the beginning of each enrichment. For the first algorithm, we start from initial offline space $V_{ms}^{n,0}$ in each time step, which means the residual-driven basis functions obtained in last enrichment are not used in current time step. On the other hand, for the second one, the residual-driven basis functions from previous time steps are kept and will be used afterwards. Using this accumulation strategy in the second algorithm, we can skip offline enrichment after a certain time period when the residual defined in (\ref{loc res}) is under a given tolerance.\\

	\begin{table}[!htbp]
		\begin{tabular}{l}
			
			\hline\hline
			\textbf{OFFLINE STAGE 1:} \\
			Construction of offline multiscale space $V_{ms}$. \\
			
			\textbf{OFFLINE STAGE 2:} On specific time step $t_n$, offline enrichment start from $V_{ms}$.\\
			Fix the number $\theta$ with $0\leq \theta\leq 1$ and error tolerance $\tau>0$.\\
			Let $V_{ms}^{n,0}=V_{ms}$.\\
			for k-th enrichment\\
			\textbf{Input}:  $V_{ms}^{n,k-1}$.\\
			1. Find the multiscale solution. Solve the reduced system (\ref{ms sol}) in $V_{ms}^{n,k-1}$.\\
			Compute the error and break the for loop if the error decreases below $\tau$.\\
			2. Compute the error indicators and select coarse neighborhoods where basis enrichment is needed.\\
			Compute total residual $R_k$ (defined in (\ref{glo res})),  the local residuals in each local domain (defined in (\ref{loc res})) . \\
			Find $r_{k_1},\ldots,r_{k_m}$, the biggest $m$ local residuals such that $r_{k_1}+\cdots+r_{k_m}\geq \theta R_k$, where $\theta$ is a given fraction. \\
			3. Solve the residual-driven basis and add them to the space.\\
			Solve (\ref{online}) in the above $k$ local domains with above residuals $r_{k_1},\cdots,r_{k_m}$ respectively to obtain $\beta_{k_1}^{n,k},\ldots,\beta_{k_m}^{n,k}$.\\
			$V_{ms}^{n, k}=V_{ms}^{n, k-1} \oplus \operatorname{span}\{\beta_{k_1}^{n,k},\ldots,\beta_{k_m}^{n,k}\}$.\\
			
			\textbf{Output}:  $V_{ms}^{n,k}$.\\
			end for\\
			\hline
		\end{tabular}
		\caption{Adaptive offline enrichment method 1.}
		\label{online algo1}
	\end{table}

	\begin{table}[!htbp]
		\centering
		\begin{tabular}{l}
			\hline\hline
			\textbf{OFFLINE STAGE 1:} \\
			Construction of offline multiscale space $V_{ms}$. \\
			
			\textbf{OFFLINE STAGE 2:} On specific time step $t_n$, offline enrichment start from $V_{ms}$.\\
			Fix the number $\theta$ with $0\leq \theta\leq 1$ and error tolerance $\tau>0$.\\
			Suppose $V_{ms}^{n-1}$ is the multiscale space in time $t_{n-1}$.\\
			Let $V_{ms}^{n,0}=V_{ms}^{n-1}$, which is the \textbf{only} difference from adaptive offline enrichment method 1.\\
			for k-th enrichment\\
			\textbf{Input}:  $V_{ms}^{n,k-1}$.\\
			1. Find the multiscale solution. Solve the reduced system (\ref{ms sol}) in $V_{ms}^{n,k-1}$.\\
			Compute the error and break the for loop if the error decreases below certain error tolerance $\tau$.\\
			2. Compute the error indicators and select coarse neighborhoods where basis enrichment is needed.\\
			Compute total residual $R_k$ (defined in (\ref{glo res})),  the local residuals in each local domain (defined in (\ref{loc res})) . \\
			Find $r_{k_1},\ldots,r_{k_m}$, the biggest $m$ local residuals such that $r_{k_1}+\cdots+r_{k_m}\geq \theta R_k$, where $\theta$ is a given fraction. \\
			3. Solve the residual-driven basis and add them to the space.\\
			Solve (\ref{online}) in the above $k$ local domains with above residuals $r_{k_1},\cdots,r_{k_m}$ respectively to obtain $\beta_{k_1}^{n,k},\ldots,\beta_{k_m}^{n,k}$.\\
			$V_{ms}^{n, k}=V_{ms}^{n, k-1} \oplus \operatorname{span}\{\beta_{k_1}^{n,k},\ldots,\beta_{k_m}^{n,k}\}$.\\
			
			\textbf{Output}:  $V_{ms}^{n,k}$.\\
			end for\\
			\hline
		\end{tabular}
		\caption{Adaptive offline enrichment method 2. }
		\label{online algo2}
	\end{table}

	\section{GMsFEM-POD method}
	In this section, we  introduce our local-global method. Different from the deterministic case in Section \ref{sec:GMsFEM}, we will consider the target random field $\kappa(\mathbf{x},\omega)$  in this section, where $\omega$ is a random variable. Before introducing the method GMsFEM-POD, we briefly review two techniques: Karhunen-Lo$\grave{\mathbf{e}}$ve expansion \cite{huang2001convergence} and proper orthogonal decomposition \cite{chatterjee2000introduction}.
	\subsection{Karhunen-Lo$\grave{\mathbf{e}}$ve expansion}\label{sect:KLE}
	Karhunen-Lo$\grave{e}$ve expansion (KLE) is an efficient method of representing random permeability fields. Let 
	$\kappa(\mathbf{x},\omega)$ be a random field. To guarantee positive permeability almost surely in concerned domain $\Omega$, we consider the natural logarithm of permeability, i.e. $Y=\ln{\kappa}$. In particular, $Y(\mathbf{x},\omega)=\bar{Y}(\mathbf{x})+Y'(\mathbf{x},\omega)$, where $\bar{Y}(\mathbf{x})=\mathrm{E}(Y(\mathbf{x},\omega))$. We further define the mean permeability field $\bar{\kappa}(\mathbf{x})$ by $\bar{\kappa}(\mathbf{x})=\exp(\bar{Y}(\mathbf{x}))$. To obtain the KLE of $Y$, we need to solve an eigenvalue problem using the covariance function $C_Y(\mathbf{x},\mathbf{z})$. In this paper, we choose a two-point exponential covariance function $C_Y$ used in \cite{he2019reduced}, that is,
	\begin{eqnarray}
		C_Y(\mathbf{x},\mathbf{z})=\sigma_{Y}^2\exp\big[-\dfrac{|x_1-z_1|^2}{\eta_1^2}-\dfrac{|x_2-z_2|^2}{\eta_2^2}\big],
		\label{KLE}
	\end{eqnarray}
	
	where $\mathbf{x}=(x_1,x_2),  \quad \mathbf{z}=(z_1,z_2)$, $\sigma_{Y}^2$ is the variance of the stochastic field $Y$. $\eta_1$ and $\eta_2$ are defined to be the correlation lengths in the first and second directions. We then solve eigenvalues and eigenvectors via the following integral equation:
	\begin{eqnarray}
		\int_{\Omega}C_Y(\mathbf{x},\mathbf{z})f(\mathbf{z})d\mathbf{z}=\lambda f(\mathbf{x}).
	\end{eqnarray}
	We then choose a set of dominating eigenvalues and corresponding eigenvectors denoted as $\{\lambda_i\}_{i=1}^{N}$ and $\{f_i\}_{i=1}^{N}$, where $N$ is the number of chosen eigenvalues.
	Subsequently, the truncated Karhunen-Lo$\grave{\mathbf{e}}$ve expansion of the log permeability can be expressed as follows,
	\begin{eqnarray}
		Y(\mathbf{x},\omega)=E[Y](\mathbf{x})+\Sigma_{i=1}^{N}\eta_i\sqrt{\lambda_i}f_i(\mathbf{x}),
	\end{eqnarray}
	where $\{\eta_i\}$ are mutually uncorrelated random variables with zero mean and unit variance. When $Y'$ is given by Gaussian process, $\{\eta_i\}$ are independent. With a fixed set of eigenvalues and eigenvectors, we can generate a set of realizations of the random permeability field $\kappa(\mathbf{x},\omega)$ by using different coefficients $\{\eta_i\}$.
	\subsection{Proper orthogonal decomposition}
	\label{sect: POD}
	In this subsection, we  introduce the proper orthogonal decomposition (POD) method. 
	
	For all $u\in H_0^1(\Omega)$, we can define a corresponding vector $\tilde{U}\in \mathrm{R}^m$ such that $u(x_i)=\tilde{U}_i$ where $\{x_i\}$ are nodal points on fine mesh $\mathcal{T}_h$.
	Besides, we define finite difference quotients $\Delta u_{H}^{j}$ by $\Delta u_{H}^{j}=\frac{u_{H}^{j}-u_{H}^{j-1}}{\Delta t}$, for $1\leq j\leq n$.
	To apply POD, we first construct a global snapshot space. Suppose we have a specific sample permeability field $\kappa(\mathbf{x},\omega_i)$, which is selected according to Theorem \ref{final_thm}. We define $\{u_{H}^{j}(\omega_i)\}\in V_{ms}$ to be the multiscale solution to $\eqref{ms sol stoc}$ associated with $\omega_i$ at time $t_j$.
	 We obtain a set of solution samples $\{u_{H}^{j}(\omega_i)\}\in V_{ms}$ and finite difference quotients $\{\Delta u_{H}^{j}(\omega_i)\}$ for $1\leq j\leq n$. We further define snapshot space $W$ as follows.
	\begin{align*}
	W=\{u_{H}^{j}(\omega_i),\Delta{u}_{H}^{j}(\omega_i)\}, \quad \text{ for }1\leq j\leq n.
	\end{align*}
	  Since $\omega_i$ is fixed, we here omit it when no ambiguity occurs. Recall $M$ is the dimension of multiscale space after offline enrichment (see Table \ref{online algo1} and \ref{online algo2}). Here we have $n\geq M$. $W$ is called solution snapshot space, which should be distinguished from the previously mentioned $V_{snap}$ spanned by local snapshots constructed in (\ref{snap_basis}). To obtain a set of POD basis from $W$, we further solve the following spectral problem:
	\begin{eqnarray}
		\sum_{i=1}^{n}\langle u_{H}^{i},\psi\rangle_{a}u_{H}^{i}+Q^2\sum_{i=1}^{n}\mathcal{A}(\Delta u_{H}^{i},\psi)\Delta u_{H}^{i}=\lambda \psi,\label{eigen_pod_con}
	\end{eqnarray}
	where $Q$ is defined in (\ref{elliptic_L2}).
	To solve (\ref{eigen_pod_con}) numerically, we need to derive the matrix form of it.
	
	Let $\{c_{H}^{i}\}\in \mathrm{R}^m$ be the vector from expanding $\{u_{H}^{i}\}\in V_{ms}$ at nodal points, i.e. $u_{H}^{i}(x_j)=c_{H}^{i}[j]$ where $\{x_j\}$ are nodal points on fine mesh $\mathcal{T}_h$ and $c_{H}^{i}[j]$ is the $i$-th entry of $c_{H}^{i}$. We can further define  $\{\Delta c_{H}^{i}\}$ by letting $\Delta c_{H}^{i}=\frac{c_{H}^{i}-c_{H}^{i-1}}{\Delta t}$, for $1\leq i\leq n$. Let $Y_1=[c_{H}^{1},\ldots,c_{H}^{n}]\in \mathrm{R}^{m*n}$,
	and $Y_2=[\Delta c_{H}^{1},\ldots,\Delta c_{H}^{n}]\in \mathrm{R}^{m*n}$ and we can have discretized version of (\ref{eigen_pod_con}) as follows,
	\begin{eqnarray}
		(Y_1Y_1^TA+Q^2Y_2Y_2^TA) \phi =\lambda \phi. \label{eigen_pod}
	\end{eqnarray}
	Comparing (\ref{eigen_pod_con}) and (\ref{eigen_pod}), we have that $\phi\in \mathrm{R}^m$ is discretized version of $\psi(\mathrm{x})$. 
	
	Since $M$ is the dimension of multiscale space, we can obtain a set of eigenpairs $\{\lambda_i,\psi_i\}_{1\leq i\leq M}$. We arrange the eigenvalues in a descending order and we have $\lambda_1\geq \lambda_2\geq \cdots \lambda_{M}$ while the corresponding eigenvectors are $\psi_1,\cdots,\psi_{M}$ respectively. 
	Without loss of generality, we assume the eigenvectors are orthonormal. In particular, $\langle \psi_i,\psi_j\rangle_{a}=\delta_{ij}$ for $1\leq i,j\leq M$. We take $l$ eigenvectors corresponding to the $l$ dominating eigenvalues to obtain the POD space, where $l<M$. We use $Q_{l,W}$ denote the POD space spanned by $l$ dominating eigenvectors resulted from W, i.e. $Q_{l,W}=\text{span}\{\psi_1,\ldots,\psi_l\}$. We further define the orthogonal projection operator $S_{\{\psi_1,\ldots,\psi_l\}}: H^1(\Omega)\to Q_{l,W}$ as follows. For each $y\in H^1(\Omega)$,
	\begin{eqnarray}
		S_{\{\psi_1,\ldots,\psi_l\}} y=\sum_{i=1}^{l}\mathcal{A}( y,\psi_i).
	\end{eqnarray}
	\begin{lemma}\label{lemma:POD}\cite{liang2002proper}
		Assume we have obtained $W$ and a set of POD basis $\{(\lambda_i,\psi_i)\}_{i=1}^{l}$ by solving (\ref{eigen_pod}), then we have $\{\psi_i\}_{i=1}^{l}$ are the solutions to the following problem
		\begin{eqnarray}
			\begin{aligned}
				&\min_{\tilde{\psi_1},\ldots,\tilde{\psi_l}}\sum_{i=1}^{n}\|u_{H}^{i}-
				S_{\{\tilde{\psi_1},\ldots,\tilde{\psi_l}\}}u_{H}^{i}\|_{a}^2+Q^2\sum_{i=1}^{n}\|\Delta u_{H}^{i}-
				S_{\{\tilde{\psi_1},\ldots,\tilde{\psi_l}\}} \Delta u_{H}^{i}\|_a^2,\\
				& s.t. \langle \tilde{\psi}_i,\tilde{\psi}_j \rangle_{a}=\delta_{ij}
			\end{aligned}
		\end{eqnarray}
		Moreover, with $S_{\{\psi_1,\ldots,\psi_l\}}$,
		\begin{eqnarray}
			\displaystyle\sum_{i=1}^{n}\|u_{H}^{i}-
			S_{\{\psi_1,\ldots,\psi_l\}} u_{H}^{i}\|_{a}^2+Q^2\sum_{i=1}^{n}\|\Delta u_{H}^{i}-
			S_{\{\psi_1,\ldots,\psi_l\}}\Delta u_{H}^{i}\|_{a}^2
			=\displaystyle\sum_{p=l+1}^{M}\lambda_p.\label{error:POD}
		\end{eqnarray}
	\end{lemma}
	Remark:
	From this lemma, we can adjust $l$ to set threshold for approximation. 
	
	\subsection{GMsFEM-POD}\label{sect:GMsFEM-POD}
	In this subsection, we introduce the local-global method: GMsFEM-POD, which is a combination of GMsFEM and POD method. The construction of global basis consists of three steps: first, apply GMsFEM to obtain multiscale basis (local basis); secondly, solve (\ref{ms sol stoc}) with different samples $\kappa(\mathbf{x},\omega)$ in multiscale space obtained in first step and use the obtained multiscale solutions to construct global snapshot space $W$; thirdly, perform POD on $W$ to obtain POD basis.  In terms of offline enrichment in the first step, we compare different offline enrichment algorithms as elaborated in Tables \ref{tab:algo1} and \ref{tab:algo2} respectively. The only difference exists in permeability field used in offline enrichment  in the first step. In the first algorithm, we perform both offline stage 1 and 2 on the mean permeability field $\bar{\kappa}(\mathbf{x})$. On the other hand, for the second algorithm, we perform offline enrichment hierarchically based on some well-chosen samples. Recall that $\bar{\kappa}(\mathbf{x})$ is the mean permeability field. Suppose we have a set of samples of permeability fields $\kappa(\mathbf{x},\omega_1),\ldots,\kappa(\mathbf{x},\omega_k)$. We perform the offline stage 1 on the $\bar{\kappa}(\mathbf{x})$. And then we perform the offline stage 2 with $\kappa(\mathbf{x},\omega_i)$ for $1\leq i\leq k$ and compute global residual $R_1,\ldots R_k$, which is defined in (\ref{glo res}). We perform the offline enrichment with the permeability field which corresponds to the largest residual norm. Using this strategy, one can make use of the information contained in the samples.
	\begin{table}[htbp!]
		\centering
		\begin{tabular}{c l}
			\hline \hline
			\textbf{Step1}& Construct multiscale basis space associated with mean permeability field $\bar{\kappa}(\mathbf{x})$. \\
			\textbf{Step2}& Solve (\ref{ms sol stoc}) associated with some samples $\kappa(\mathbf{x},\omega)$ using multiscale basis in step 1\\
			& and obtain solution snapshot space $W$.\\
			\textbf{Step3}&  Apply POD on $W$ to obtain POD basis and solve (\ref{snap sol stoc}) for each $\omega$. \\
			\hline
		\end{tabular}
		\caption{GMsFEM-POD method 1.}
		\label{tab:algo1}
	\end{table}
	\begin{table}[htbp!]
		\centering
		\begin{tabular}{c l}
			\hline\hline
			\textbf{Step1}& \textbf{Offline stage 1} \\
			& Apply GMsFEM on mean permeability $\bar{\kappa}(\mathbf{x})$ to obtain multiscale space $\Phi_{0}$.  \\
			& \textbf{Offline stage 2} \\
			& Use KLE to obtain N realizations of permeability fields $\{\kappa_i\}_{ 1\leq i \leq N}$. \\
			&Let $V_{ms}^{0}=V_{ms}$.\\
			&for k-th enrichment\\
			&\textbf{Input}:  $V_{ms}^{k-1}$.\\
			& \quad \quad  For each $1\leq i \leq N$, solve (\ref{ms sol}) associated with $\kappa_i$ and compute $\|R_i\|$ with (\ref{glo res}). \\
			& \quad \quad Let $\|R_{p_k}\|=\max\{\|R_{p_k}\|,\ldots,\|R_{N}\|\}.$\\
			& \quad \quad Break the for loop if $\|R_{p_k}\|$ decreases below a given tolerance.\\
			&\quad \quad Perform online enrichment with the permeability field $\kappa_{p_k}$ to obtain $V_{ms}^{k}$.\\
			
			&\textbf{Output:} $V_{ms}^{k}$.\\
			\textbf{Step2 \& Step3} &Same as GMsFEM-POD method 1.\\
			\hline
		\end{tabular}
		\caption{GMsFEM-POD method 2.}
		\label{tab:algo2}
	\end{table}
	\section{Error analysis}
	In this section, we will present error analysis. First of all, we will introduce some notations. 
	We will consider the random case and we use $\omega$ to denote the random variable in $\Omega_r$ subject to a specific distribution. 
	We define the $V_{ms,\omega}$ to be the multiscale space constructed using $\kappa(\mathbf{x},\omega)$. In particular, we solve \eqref{snap_basis}-\eqref{POU} using $\kappa(\mathbf{x},\omega)$. $Q_{l,W}$ is defined to be the POD space based on some solution space $W$ as discussed in Section \ref{sect: POD}. Let $u_h^n(\omega) \in V_h$, $u_H^n(\omega) \in V_{ms,\omega}$ and $p_{l}^{n}(\omega)\in Q_{l,W}$ be the solutions to (\ref{fine sol stoc}), (\ref{ms sol stoc}) and (\ref{snap sol stoc}) respectively associated with permeability field $\kappa(\mathbf{x},\omega)$ at $t=t_n$. Let $\kappa_{\text{min}}(\omega)=\min_{\mathbf{x}}{\kappa(\mathbf{x},\omega)}$ for each $\omega\in \Omega_r$. Recall that $f^n=f(t_n,\mathbf{x})$.
	
	For each $\omega \in \Omega_r$, we have
	\begin{eqnarray}
		\begin{aligned}
			\frac{1}{\Delta t}\left\langle u_h^n(\omega), v\right\rangle+\left\langle\kappa(\mathbf{x},\omega) \nabla u_h^n(\omega), \nabla v\right\rangle&=\left\langle \frac{1}{\Delta t} u_{h}^{n-1}(\omega)+f^n, v\right\rangle, \\
			\langle u_{h}^{0}(\omega),v\rangle&=\langle g,v\rangle,\quad \forall v \in V_h.
		\end{aligned}
		\label{fine sol stoc}
	\end{eqnarray}
	\begin{eqnarray}
		\begin{aligned}
			\frac{1}{\Delta t}\left\langle u_H^n(\omega), v\right\rangle+\left\langle\kappa(\mathbf{x},\omega)\nabla u_H^n(\omega), \nabla v\right\rangle&=\left\langle \frac{1}{\Delta t} u_{H}^{n-1}(\omega)+f^n, v\right\rangle, \quad \\
			\langle u_{H}^{0}(\omega),v\rangle&=\langle g,v\rangle, \quad\forall v \in V_{ms,\omega}.\label{ms sol stoc}
		\end{aligned}
	\end{eqnarray}
	
	\begin{eqnarray}
		\begin{aligned}
			\frac{1}{\Delta t}\left\langle p_{l}^{n}(\omega), v\right\rangle+\left\langle\kappa(\mathbf{x},\omega)\nabla p_{l}^{n}(\omega), \nabla v\right\rangle&=\left\langle \frac{1}{\Delta t} p_{l}^{n-1}(\omega)+f^n, v\right\rangle, \quad \\
			\langle p_{l}^{0}(\omega),v\rangle&=\langle g,v\rangle, \quad\forall v \in Q_{l,W}.\label{snap sol stoc}
		\end{aligned}
	\end{eqnarray}
	Our goal here is to estimate the difference between the $p_{l}^{n}(\omega)$ and our reference solution $u_h^n(\omega)$ under some specific norm and we split the error into four parts as follows:
	\begin{eqnarray}
		\begin{aligned}
			u_h^n(\omega)- p_{l}^{n}(\omega)&=u_h^n(\omega)-u_h^n(\omega_i)+u_h^n(\omega_i)-u_H^n(\omega_i)\\
			&+u_H^n(\omega_i)-p_{l}^{n}(\omega_i)+p_{l}^{n}(\omega_i)-p_{l}^{n}(\omega)\\
			&:=e_{1}^{n}+e_{2}^{n}+e_{3}^{n}+e_{4}^{n},\label{split error}
		\end{aligned}
	\end{eqnarray}
	where $\omega_i$ is a sample in $\Omega_r$ used to construct the POD basis (see Section \ref{sect: POD}).
	
	In the following part, unless specified, the norm is calculated in domain $\Omega$. Also, we  suppress $\Omega$ when no ambiguity occurs. To evaluate the above four error terms, we make several assumptions below.

	\begin{assumption}
		Suppose the random permeability field $\kappa(\mathbf{x},\omega)$ satisfies the following property:
		given $\omega\in \Omega_r$ and $\epsilon_1$, there exists an integer $N_{\epsilon_1}$ and a set of $\{\kappa(\mathbf{x},\omega_i)\}$, $i=1,\ldots,N_{\epsilon_1}$ such that
		\begin{eqnarray}
			\mathrm{E}[\inf_{1\leq i\leq N_{\epsilon_{1}}}\|\kappa(\mathbf{x},\omega)-\kappa(\mathbf{x},\omega_i)\|_{L^{\infty}(\Omega)}]\leq \epsilon_1.
		\end{eqnarray}\label{ass1}
	\end{assumption}
	This assumption characterizes a "density" property of the snapshot space $\{\kappa(\mathbf{x},\omega)\}$, which is necessary in the construction of snapshot space for solutions. The verification of this assumption can be seen in (29) in \cite{li2020data}.
	
	\begin{assumption}\label{ass3}
		We assume that there is a positive constant $C_1$ such that the force function $f$ in {\rm (\ref{eqn:model})} satisfies
	\begin{align}
		\Delta t\sum_{j=1}^{n}\|f^j\|_{L^2}\leq C_1,
	\end{align}
	\end{assumption}
   where $f^j=f(t_j,\mathbf{x})$, for $1\leq j\leq n$.
	\begin{lemma}\label{stability}
		Suppose $u_h^i(\omega)$ to be the solution to \eqref{fine sol stoc} at $t=t_i$ for $i=1,\ldots, n$. Then, for each $\omega\in \Omega_r$, we have 
		\begin{eqnarray}
	\|u_h^n(\omega)\|_{L^2}+\sqrt{\frac{\kappa_{min}(\omega)}{Q}}\Delta t\sum_{j=1}^{n}|u_h^j(\omega)|_{H_0^1}\leq \|u_h^0\|_{L^2}+\Delta t
\sum_{j=1}^{n}\|f^j\|_{L^2}.\label{ass2_result}
	\end{eqnarray}
where $Q$ is defined in \eqref{elliptic_L2}.
	\end{lemma}
	\begin{proof}
	With (\ref{fine sol stoc}), we have
	\begin{eqnarray}
 \langle u_h^n(\omega),v\rangle +\Delta t \mathcal{A}(u_h^n(\omega),v)=\langle u_h^{n-1}(\omega),v\rangle+\Delta t\langle f,v\rangle, \quad\forall v\in V_{f}.
	\end{eqnarray}
Take $v=u_h^n(\omega)$ and we have
	\begin{eqnarray*}
	\langle u_h^n(\omega),u_h^n(\omega)\rangle +\Delta t \mathcal{A}(u_h^n(\omega),u_h^n(\omega))=\langle u_h^{n-1}(\omega),u_h^n(\omega)\rangle+\Delta t\langle f,u_h^n(\omega)\rangle, \quad \forall v\in V_{f}.
\end{eqnarray*}
Applying Cauchy-Schwartz inequality, we have
\begin{align*}
	\|u_h^n(\omega)\|_{L^2}+\Delta t \frac{\|u_h^n(\omega)\|_a^2}{\|u_h^n(\omega)\|_{L^2}}\leq \|u_h^{n-1}(\omega)\|_{L^2}+\Delta t \|f^n\|_{L^2}.
	\end{align*}
Using \eqref{elliptic_L2}, we have ,
\begin{align}
	\|u_h^j(\omega)\|_{L^2}+\sqrt{\frac{\kappa_{min}(\omega)}{Q}}\Delta t |u_h^j(\omega)|_{H_0^1}\leq \|u_h^{j-1}(\omega)\|_{L^2}+\Delta t \|f^j\|_{L^2},\quad 
	\forall j=1,\ldots, n. \label{add_lemma}
\end{align}
Adding \eqref{add_lemma} for $j=1,\ldots, n$, we have \eqref{ass2_result}.
	\end{proof}
	Next  we present the error estimates for $e_1$ and $e_4$. We remark that $u$ depends on $\omega$. For simplicity of notation, we let $\kappa(\omega)=\kappa(\mathbf{x},\omega)$ for each $\omega\in \Omega_r$.
	
	\begin{theorem}\label{e_1}
		Under the above Assumptions, we can have the following error estimates, where $\omega \in \Omega_r$ and $\omega_i$ is a sample used in the snapshot space in the Section \ref{sect: POD}.
		\begin{align}
			\|e_{1}^{n}\|_{L^2}^2&\leq \frac{C}{\kappa_{min}(\omega)\kappa_{min}(\omega_i)}
			\|\kappa(\omega)-\kappa(\omega_i)\|_{L^{\infty}}^2,\label{e_1_bd}\\
			\|e_{4}^{n}\|_{L^2}^2&\leq \frac{C}{\kappa_{min}(\omega)\kappa_{min}(\omega_i)}
			\|\kappa(\omega)-\kappa(\omega_i)\|_{L^{\infty}}^2,\label{e_4_bd}
		\end{align}
		where $e_{1}^{n}=u_h^n(\omega)-u_h^n(\omega_i)$ and $e_{4}^{n}=p_l^n(\omega_i)-p_l^n(\omega)$ as in {\rm(\ref{split error})}.
	\end{theorem}
	\begin{proof}
		From (\ref{fine sol stoc}), for all $v\in V_h$, we have
		\begin{eqnarray}
			\frac{1}{\Delta t}\left\langle u_h^n(\omega), v\right\rangle+\left\langle\kappa(\omega) \nabla u_h^n(\omega), \nabla v\right\rangle&=\left\langle \dfrac{1}{\Delta t} u_{h}^{n-1}(\omega)+f^{n}, v\right\rangle,\label{fine1}\\
			\frac{1}{\Delta t}\left\langle u_h^n(\omega_i), v\right\rangle+\left\langle\kappa(\omega_i) \nabla u_h^n(\omega_i), \nabla v\right\rangle&=\left\langle \dfrac{1}{\Delta t} u_{h}^{n-1}(\omega_i)+f^{n}, v\right\rangle .\label{fine2}
		\end{eqnarray}
		Substract (\ref{fine2}) from (\ref{fine1}) and we have
		\begin{eqnarray*}
			\langle e_{1}^{n},v\rangle-\langle e_{1}^{n-1},v\rangle+\Delta t\langle \kappa(\omega)\nabla e_{1}^{n},\nabla v\rangle =-\Delta t
			\langle(\kappa(\omega)-\kappa(\omega_i))\nabla u_h^n(\omega_i),\nabla v\rangle.
		\end{eqnarray*}
		Take $v=e_{1}^{n}$, we have
		\begin{eqnarray*}
			\Delta t\langle\kappa(\omega)\nabla e_{1}^{n},\nabla e_{1}^{n}\rangle+\langle e_{1}^{n},e_{1}^{n}\rangle=\langle e_{1}^{n-1},e_{1}^{n}\rangle-\Delta t\langle(\kappa(\omega)-\kappa(\omega_i))\nabla u_h^n(\omega_i),\nabla e_{1}^{n}\rangle.
		\end{eqnarray*}
		We use the Cauchy-Schwartz inequality to obtain
		\begin{align*}
			&\kappa_{min}(\omega)\Delta t|e_{1}^{n}|_{H_0^1}^2+\|e_1^{n}\|_{L^2}^2\leq\frac{1}{2}\|e_1^{n}\|_{L^2}^2+\frac{1}{2}\|e_1^{n-1}\|_{L^2}^2+\Delta t\|(\kappa(\omega)-\kappa(\omega_i))\|_{L^{\infty}}\|u_h^n(\omega_i)\|_{H_0^1}\|e_1^n\|_{H_0^1}\\
			&\leq\frac{1}{2}\|e_1^{n}\|_{L^2}^2+\frac{1}{2}\|e_1^{n-1}\|_{L^2}^2+\frac{\Delta t}{4\kappa_{min}(\omega)}\|(\kappa(\omega)-\kappa(\omega_i))\|_{L^{\infty}}^2|u_h^n(\omega_i)|_{H_0^1}^2+\kappa_{min}(\omega)\Delta t|e_1^n|_{H_0^1}^2.
		\end{align*}
		Hence, for $1\leq j\leq  n$, we have
		\begin{align}\label{e_1_ele}
			\|e_1^{j}\|_{L^2}^2\leq \|e_1^{j-1}\|_{L^2}^2+\frac{\Delta t}{2\kappa_{min}(\omega)}\|(\kappa(\omega)-\kappa(\omega_i))\|_{L^{\infty}}^2|u_h^j(\omega_i)|_{H_0^1}^2.
		\end{align}
		Adding (\ref{e_1_ele}) for $1\leq j\leq  n$ and using Lemma \ref{stability}, we have 
		\begin{align*}
			\|e_1^{n}\|_{L^2}^2\leq \frac{\Delta t}{2\kappa_{min}(\omega)}\|(\kappa(\omega)-\kappa(\omega_i))\|_{L^{\infty}}^2\sum_{j=1}^{n}|u_h^j(\omega_i)|_{H_0^1}^2.
		\end{align*}
		Here we use the fact that $e_{1}^{0}=0$. Apply Lemma \ref{stability} and we have $\sum_{j=1}^{n}|u_h^j(\omega)|_{H_0^1}^2\leq \frac{C}{\kappa_{min}(\omega_i)}$. Hence (\ref{e_1_bd}) holds. We use the similar technique as above to obtain (\ref{e_4_bd}).\\
		
	\end{proof}
	
	In the remaining part, we focus on estimating the errors $e_2$ and $e_3$ in (\ref{split error}). Since we consider the error with sample $\omega_i$ which is a deterministic case, we use $\kappa:=\kappa(\mathbf{x},\omega_i)$ when no ambiguity occurs.
	\begin{definition}\label{ritz}
		We define $R_H: H_{0}^{1}(\Omega)\rightarrow V_{ms}$ to be the elliptic projection 
		with respect to the energy inner product, more explicitly
		$$ \langle\kappa \nabla(R_H u-u),\nabla v\rangle=0, \quad\forall v\in V_{ms}. $$
		
	\end{definition}
	\begin{lemma}\label{refletion}
		\cite{li2019convergence} Let $u\in H_0^1(\Omega)$ be solution to the elliptic equation
		\begin{eqnarray}
			\left\{\begin{array}
			{rlrl}{-d i v(\kappa \nabla u)} & {=f} & {} & {\text {\rm in } \Omega}, \\
				{u} & {=0} & {} & {\text {\rm on } \partial \Omega.}	\end{array}\right.\label{elliptic problem}
		\end{eqnarray}
	 Define $u_h$ to be the fine-scale solution to \eqref{elliptic problem}. Specifically, 
	 \begin{align}
	 	\mathcal{A}(u_h,v)=\langle f,v\rangle\quad  \forall v\in V_h.
	 \end{align}
		If $R_{H}u_h$ is the elliptic projection of $u_h$ in the subspace $V_{ms}$, we have 
			\begin{align}
		\|u_h-R_H u_h\|_a&\leq \gamma(H)\|\kappa^{-\frac{1}{2}}f\|_{L^2}, \label{energ_ell_pro}\\ 
		\|u_h-R_{H}u_h\|_{L^2}&\leq (\gamma(H))^2\kappa_{min}({\omega_i})^{-\frac{1}{2}}\|\kappa^{-\frac{1}{2}}f\|_{L^2}\label{l2_ell_pro}.
			\end{align}
		where $\kappa_{min}(\omega_i)=\min_{\mathbf{x}}\{\kappa(\mathbf{x},\omega_i)\}$.
		Moreover $$\gamma(H)=\sqrt{2C_{ov}}H\max_{i=1,\ldots,N_{in}}\{C_0HC_{poin}(D_i)+\sqrt{C_{poin}(D_i)}\}+
		C_1\sqrt{20C_{ov}}(H^2\Lambda)^{-\frac{1}{2}}.$$ The concerned constants are defined as follows:
		\begin{align}
			C_{ov}&:=\max_{K\in\mathcal{T}^{H}}\#\{O_i:K\subset D_i\text{ for } i=1,\ldots,N_{in}\}.\\
			C_{poin}(D_i)&:=H^{-2}\max_{\eta\in H_{\kappa}(D_i)}\frac{\int_{D_i}\kappa\eta^2}{\int_{D_i}\kappa|\nabla\omega|^2}.\\
			C_{poin}(\Omega)&:=\text{diam}(\Omega)^{-2}\max_{\eta\in H_{\kappa}(\Omega)}\frac{\int_{\Omega}\kappa\eta^2}{\int_{\Omega}\kappa|\nabla\omega|^2}.
		\end{align}
	  $\Lambda$ is defined in \eqref{tot_eigen_val}. The constants $C_{poin}(D_i)$ and $C_{poin}(\Omega)$ are independent of $\kappa$. $C_0$ depends on $\Omega$, the size and shape of subsets $P_j$ for $j=1,\ldots,s$, the space dimension $d$ and the coefficient $\kappa$ but it is independent of the distances between the inclusions $P_k$ and $P_j$ for $k,j=1,\ldots,s$. However, the precise dependence of the constant $C_0$ on $\kappa$ is still unknown. 
	  The constant $C_1$ is given by $C_1:=H\max_{i=1,\ldots, N_{in}}\{\sqrt{C_{poin}(D_i)}\}+2diam(\Omega)\sqrt{C_{poin}(\Omega)}$.  
	\end{lemma}
\begin{proof}

Using Lemma 4.6 and Lemma 4.10 in \cite{li2019convergence}, we can get \eqref{energ_ell_pro}. 
	
Now we give estimation for $\|u-R_{H}u\|_{L^2}$ by standard dual argument.
	Suppose $w_h\in V_h$ is the solution to the following equation: 
	\begin{align}
		\mathcal{A}(w_h,v)=\langle u_h-R_{H}u_h,v\rangle, \quad\forall v\in V_{h}. \label{snap_eliptic_pro}
	\end{align}
    $R_{H}w_h$ is the elliptic projection of $w_h$. By the definition, we have 
	\begin{align}
	\mathcal{A}(R_{H}w_h,v)=\langle u_h-R_{H}u_h,v\rangle, \quad\forall v\in V_{ms}. \label{off_eliptic_pro}
     \end{align}    
Combing with \eqref{energ_ell_pro}, we have
 \begin{align}
 	\|u_h-R_H u_h\|_{L^2}^2&=\mathcal{A}(w_h,u_h-R_H u_h)=\mathcal{A}(w_h-R_H w_h,u_h-R_H u_h)\\
 	&\leq \|w_h-R_H w_h\|_a\|u_h-R_H u_h\|_a\leq(\gamma(H))^2\|\kappa^{-\frac{1}{2}}f\|_{L^2}\|\kappa^{-\frac{1}{2}}(u_h-R_H u_h)\|_{L^2}\\
 	&\leq (\gamma(H))^2\kappa_{min}(\omega_i)^{-\frac{1}{2}}\|\kappa^{-\frac{1}{2}}f\|_{L^2}\|u_h-R_H u_h\|_{L^2}.
 \end{align}
\end{proof}

	
	Before we estimate $e_2$, we first define 
	\begin{align}
		\partial_t u_h^j=\frac{u_h^{j}-u_h^{j-1}}{\Delta t}, \quad 1\leq j\leq n.\label{diff_uh}
	\end{align}
	\begin{theorem}\label{e_2}
		Let $u_h^n$ and $u_H^n$ be the solution at $t=t_n$ of equations {\rm \eqref{fine sol}} and {\rm\eqref{ms sol}}, respectively. Then, $e_{2}^{n}$ satisfies 
		\begin{eqnarray}\label{e2_1}
			\left\|e_{2}^{n}\right\|_{L^2}^2\leq C\left(\|u_h^0-u_H^0\|_{L^2}^2+\|u_h^0-R_Hu_h^0\|_{L^2}^2+\sum_{i=1}^{n}(\gamma(H))^4\kappa_{min}({\omega_i})^{-1}\|\kappa^{-\frac{1}{2}}(f^{i}-\partial_t u_h^i)\|_{L^2}^2\right). 
		\end{eqnarray}
		
	\end{theorem}
	\begin{proof}
		First, we split the error into two parts.
		\begin{eqnarray*}
			\begin{aligned}
				\left\|u_H^n-u_h^{n}\right\|_{L^2} & \leq\left\|\theta^{n}\right\|_{L^2}+\left\|\rho^{n}\right\|_{L^2}. \\ 
				\theta^{n}(\mathbf{x}) &=u_H^n-R_{H}u_h^{n}.\\ 
				\rho^{n}(\mathbf{x}) &=R_{H} u_h^{n}-u_h^{n}.
			\end{aligned} 
		\end{eqnarray*}
		We estimate the component $\rho^n$ and $\theta^n$ respectively.
	
		By Lemma \ref{refletion}
		\begin{eqnarray}
			\|\rho^n\|_{L^2}\leq (\gamma(H))^2\kappa_{min}({\omega_i})^{-\frac{1}{2}}\|\kappa^{-\frac{1}{2}}(f^{n}-\partial_t u_h^n)\|_{L^2}.\label{add1}
		\end{eqnarray}
		
		Next, we estimate $\|\theta^n\|_{L^2}$.\\
		We have the following equation for $v\in V_{ms}$.\\
		\begin{eqnarray}
			\begin{aligned}
				&\left\langle \dfrac{\theta^n-\theta^{n-1}}{\Delta t},v\right\rangle+\mathcal{A}\left(\theta^n, v\right)=\left\langle f^n,v\right\rangle-\left\langle \dfrac{R_{H}u_h^n-R_{H}u_h^{n-1}}{\Delta t},v\right\rangle-\mathcal{A}\left(R_{H}u_h^n, v\right)\\
				&=\left\langle f^n,v\right\rangle-\left\langle \dfrac{R_{H}u_h^n-R_{H}u_h^{n-1}}{\Delta t},v\right\rangle-\mathcal{A}\left(u_h^n, v\right)\\
				&=\left\langle \partial_t u_h^n-\dfrac{R_{H}u_h^n-R_{H}u_h^{n-1}}{\Delta t},v\right\rangle,
			\end{aligned}
		\end{eqnarray}
		which is equivalent to
		\begin{eqnarray}
			\left\langle\theta^n,v\right\rangle+\Delta t \mathcal{A}\left(\theta^n, v\right)=
			\left\langle\theta^{n-1},v\right\rangle+\Delta t\left\langle \partial_t u_h^n,v\right\rangle-\left\langle R_{H}u_h^n-R_{H}u_h^{n-1},v\right\rangle.
			\label{es_1}
		\end{eqnarray}
		Taking $v=\theta^n$, and applying Lemma \ref{refletion}, Young's inequality, and triangle's inequality, we obtain
		\begin{align*}
			\left\|\theta^{n}\right\|_{L^2}^{2}+\Delta t \mathcal{A}\left(\theta^n,\theta^n\right)
			\leq\left\|\theta^{n-1}\right\|_{L^2}\left\|\theta^{n}\right\|_{L^2}
			+\Delta t\|\partial_t u_h^n-\dfrac{R_{H}u_h^n-R_{H}u_h^{n-1}}{\Delta t}\|_{L^2}\left\|\theta^n\right\|_{L^2}.
		\end{align*}
		$\forall 1\leq j\leq n$, we let $z^j=\dfrac{R_{H}u_h^j-R_{H}u_h^{j-1}}{\Delta t}-\partial_t u_h^j$, then
		\begin{eqnarray}
			\label{e_2_ele}
			\|\theta^j\|_{L^2}\leq \|\theta^{j-1}\|_{L^2}+\Delta t \|z^j\|_{L^2}.
		\end{eqnarray}
		Adding (\ref{e_2_ele}) for $1\leq j\leq n$, we have
		\begin{eqnarray}\label{theta_1_sum}
			\|\theta^n\|_{L^2}\leq \|\theta^{0}\|_{L^2}+\Delta t\sum_{j=1}^{n}\|z^j\|_{L^2}.
		\end{eqnarray}
        Since $\theta_0=u_H^0-R_H u_h^0$, we have 
        \begin{align}
        	\|\theta_0\|_{L^2}\leq \|u_h^0-u_H^0\|_{L^2}+\|u_h^0-R_Hu_h^0\|_{L^2}
        \end{align}
	
		Because
		\begin{eqnarray*}
			z^i= \dfrac{R_{H}u_h^i-R_{H}u_h^{i-1}}{\Delta t}-\dfrac{u_h^i-u_h^{i-1}}{\Delta t}.
		\end{eqnarray*}
		Using Lemma \ref{refletion}, we have 
		\begin{eqnarray}
			\Delta t\sum_{i=1}^{n}\|z^i\|_{L^2}\leq \sum_{i=1}^{n}(\gamma(H))^2\kappa_{min}({\omega_i})^{-\frac{1}{2}}\|\kappa^{-\frac{1}{2}}(f^{i}-\partial_t u^i)\|_{L^2}.\label{add5}
		\end{eqnarray}
		Combining (\ref{theta_1_sum}) to (\ref{add5}), we can obtain 
		\begin{eqnarray}
			\begin{aligned}
				\|\theta^n\|_{L^2}^2\leq C\left(\sum_{i=1}^{n}(\gamma(H))^4\kappa_{min}({\omega_i})^{-1}\|\kappa^{-\frac{1}{2}}(f^{i}-\partial_t u^i)\|_{L^2}^2+\|u_h^0-u_H^0\|_{L^2}^2+\|u_h^0-R_Hu_h^0\|_{L^2}^2\right).\label{theta_e2}
			\end{aligned}
		\end{eqnarray}
		From (\ref{add1}) and (\ref{theta_e2}), we have (\ref{e2_1}).
		
	\end{proof}
	For now, we have considered $e_1$, $e_2$ and $e_4$. We then estimate $e_3$. Suppose we have solved a set of POD basis $\{(\lambda_i,\psi_i)\}_{i=1}^{l}$ by solving (\ref{eigen_pod}). For simplicity of notation, we here let $S^{l}=S_{\{\psi_1,\ldots,\psi_l\}}$.
	\begin{theorem}\label{e_3}
		Let $\{u_{H}^{i}\}$ and $\{p_{l}^{i}\}$ for $1\leq i\leq n$ be solutions of (\rm\ref{ms sol stoc}) and (\rm\ref{snap sol stoc}) respectively, then we have the following estimations,
		\begin{eqnarray}
			\|e_{3}^{n}\|_{L^2}^2\leq
			2(\Delta t+1)\sum_{p=l+1}^{M}\lambda_p,
		\end{eqnarray}
		where $\{\lambda_i\}_{i=1}^M$ are eigenvalues solved in (\ref{eigen_pod_con}).
	\end{theorem}
	\begin{proof}
		As previously, we split $e_3$ into two parts.
		\begin{eqnarray}
			u_H^n-p_{l}^{n}=d_l^{n}+q_l^{n},\label{split_e3}
		\end{eqnarray}
		where $d_l^n=u_H^n-S^lu_H^n$ and $q_l^n=S^lu_H^n-p_{l}^{n}$.
		From (\ref{error:POD}), we have
		\begin{eqnarray}
			\displaystyle\sum_{i=1}^{n}\|d_l^n\|_a^2=\sum_{i=1}^{n}\|u_{H}^{i}-
			S^l u_{H}^{i}\|_{a}^2
			\leq\displaystyle\sum_{p=l+1}^{M}\lambda_p.
		\end{eqnarray}
		From (\ref{ms sol stoc}) and (\ref{snap sol stoc}), we have $\forall v\in Q_{l,W},$ and $0\leq j\leq n-1$,
		\begin{eqnarray}
			\displaystyle\dfrac{\left\langle u_{H}^{j+1}-u_{H}^{j}, v\right\rangle}{\Delta t}+\left\langle\kappa \nabla u_{H}^{j+1}, \nabla v\right\rangle=\left\langle f, v\right\rangle, \label{ms_e3}
		\end{eqnarray}
		\begin{eqnarray}
			\displaystyle\dfrac{\left\langle p_{l}^{j+1}-p_{l}^{j}, v\right\rangle}{\Delta t}+\left\langle\kappa \nabla p_{l}^{j+1}, \nabla v\right\rangle=\left\langle f, v\right\rangle. \label{POD_e3}
		\end{eqnarray}
		Similar as \eqref{diff_uh}, we define 
		\begin{align*}
			\partial_t q_l^j=\dfrac{q_l^{j}-q_l^{j-1}}{\Delta t}, \quad 1\leq j \leq n.
		\end{align*}
		Then,
		$$\partial_t q_l^j=\dfrac{S^lu_{H}^{j}-S^lu_{H}^{j-1}}{\Delta t}-\dfrac{u_{H}^{j}-u_{H}^{j-1}}{\Delta t}+\dfrac{u_{H}^{j}-u_{H}^{j-1}}{\Delta t}-\dfrac{p_{l}^{j}-p_{l}^{j-1}}{\Delta t}.$$
		Hence we have $\forall v\in Q_{l,W}$,
		\begin{eqnarray}
			\begin{aligned}
				\langle\partial_t q_l^j,v\rangle&=\langle\dfrac{S^l u_{H}^{j}-S^l u_{H}^{j-1}}{\Delta t},v\rangle-\langle\dfrac{u_{H}^{j}- u_{H}^{j-1}}{\Delta t},v\rangle\\
				&+\langle\dfrac{u_{H}^{j}- u_{H}^{j-1}}{\Delta t},v\rangle-\langle f,v\rangle+\langle\kappa \nabla p_{l}^{j},\nabla v\rangle\\
				&=\langle\dfrac{S^l u_{H}^{j}-S^l u_{H}^{j-1}}{\Delta t},v\rangle-\langle\dfrac{u_{H}^{j}- u_{H}^{j-1}}{\Delta t},v\rangle
				+\langle\kappa \nabla (p_{l}^{j}-u_{H}^{j}),\nabla v\rangle.       
			\end{aligned}
		\end{eqnarray}
		Let $z_{l}^{j}=\langle\dfrac{S^l u_{H}^{j}-S^l u_{H}^{j-1}}{\Delta t},v\rangle-\langle\dfrac{u_{H}^{j}- u_{H}^{j-1}}{\Delta t},v\rangle$
		and we have 
		\begin{eqnarray}
			\sum_{j=1}^n\|d_l^j\|_a^2+Q^2\sum_{j=1}^{n}\|z_l^j\|_a^2=\sum_{p=l+1}^{M}\lambda_p.\label{eigen_e3}
		\end{eqnarray}
		Take $v=q_l^{j}\in Q_{l,W}$ and we have
		\begin{eqnarray}
			\langle\partial_t q_l^j,q_l^{j}\rangle+\mathcal{A}(q_l^{j},q_l^{j}) =\langle z_{l}^{j},q_l^{j}\rangle-\mathcal{A}(d_l^{j},q_l^{j}).
		\end{eqnarray}
		Applying Cauchy-Schwartz inequality, we have
		\begin{eqnarray}\label{q_ele}
			\|q_l^{j}\|_{L^2}^2+2\Delta t\|q_l^{j}\|_a^2\leq \|q_l^{j-1}\|_{L^2}^2+\Delta t \left (\|d_l^{j}\|_a^2+Q\|z_l^{j}\|_{L^2}^2+\|q_l^{j}\|_a^2+\frac{1}{Q}\|q_l^{j}\|_{L^2}^2\right),
		\end{eqnarray}
		where $Q$ is defined in \eqref{elliptic_L2}.
		Hence,
		\begin{eqnarray}\label{e_3_ele}
			\|q_l^{j}\|_{L^2}^2\leq \|q_l^{j-1}\|_{L^2}^2+\Delta t \left (\|d_l^{j}\|_a^2+Q^2\|z_l^{j}\|_{a}^2\right).
		\end{eqnarray}
		Adding \eqref{e_3_ele} for $1\leq j\leq n$, we have
		\begin{eqnarray*}
			\|q_l^{n}\|_{L^2}^2\leq \|q_l^{0}\|_{L^2}^2+\Delta t \sum_{j=1}^{n}\left (\|d_l^{j}\|_a^2+Q^2\|z_l^{j}\|_{a}^2\right).
		\end{eqnarray*}
		By (\ref{elliptic_L2}) and (\ref{eigen_e3}), one can derive
		\begin{eqnarray}
			\|q_l^{n}\|_{L^2}^2\leq
			\|q_l^{0}\|_{L^2}^2+\Delta t\sum_{p=l+1}^{M}\lambda_p.
		\end{eqnarray}
		Then we can arrive
		\begin{eqnarray}
			\|e_{3}^{n}\|_{L^2}^2\leq
			2(\Delta t+1)\sum_{p=l+1}^{M}\lambda_p,
		\end{eqnarray}
		since we can choose $q_l^{0}=0$.
	\end{proof}
	\begin{theorem} \label{final_thm}
		Almost surely for all $\omega\in \Omega_r$,
		suppose $u_h^n(\omega)$ and $p_{l}^{n}(\omega)$ are solutions to {\rm(\ref{fine sol stoc})} and {\rm (\ref{snap sol stoc})} respectively, we have for all $\epsilon>0$, there exists a corresponding $\omega_{p}$ such that
		\begin{eqnarray}
			\begin{aligned}
				\|u_h^n(\omega)-&p_{l}^{n}(\omega)\|_{L^2}^2\leq C\left(\|u_h^0-u_H^0\|_{L^2}^2+\|u_h^0-R_Hu_h^0\|_{L^2}^2+\sum_{i=1}^{n}(\gamma(H))^4\kappa_{min}({\omega_p})^{-1}\|\kappa^{-\frac{1}{2}}(f^{i}-\partial_t u_h^i)\|_{L^2}^2\right)\\
				&+C(\Delta t+1)\sum_{p=l+1}^{M}\lambda_p+\dfrac{C\epsilon^2}{\kappa_{min}(\omega)\kappa_{min}(\omega_p)}.\label{final_1}  \end{aligned}
		\end{eqnarray}
		
	\end{theorem}
	\begin{proof}
		Based on Assumption \ref{ass1} and Theorem \ref{e_1}, almost surely for all $\omega\in \Omega_r$ , for every positive $\epsilon$, one can find a set of snapshot $\{\kappa(\mathbf{x},\omega_i)\}_{1\leq i\leq N_{\epsilon} }$ such that the following is satisfied. 
		\begin{eqnarray}
			\inf_{1\leq i\leq N_{\epsilon}}\|\kappa(\mathbf{x},\omega)-\kappa(\mathbf{x},\omega_i)\|_{L^{\infty}(\Omega)}\leq \epsilon.
		\end{eqnarray}
		Hence there exits a corresponding $\omega_{p}$ where $1\leq p\leq N_{\epsilon}$ such that
		\begin{align*}
		\|\kappa(\mathbf{x},\omega)-\kappa(\mathbf{x},\omega_p)\|_{L^{\infty}(\Omega)}\leq \epsilon.
		\end{align*}
	 Combining with Theorem \ref{e_1}, we have
		\begin{eqnarray}
			\begin{aligned}
				\|u_h^n(\omega)-u_h^n(\omega_p)\|_{L^2}^2 
				\leq \dfrac{C}{\kappa_{min}(\omega)\kappa_{min}(\omega_p)}\|\kappa(\omega)-\kappa(\omega_p)\|_{L^{\infty}}^2\leq \dfrac{C\epsilon^2}{\kappa_{min}(\omega)\kappa_{min}(\omega_p)},\label{split1}
			\end{aligned}
		\end{eqnarray}
		where we apply Poincar$\acute{e}$ inequality. 
		
		Similarly as in deriving (\ref{split1}), we have
		\begin{eqnarray}
			\|p_l^n(\mathbf{x},\omega)-p_l^n(\mathbf{x},\omega_p)\|_{L^2}^2\leq  \dfrac{C\epsilon^2}{\kappa_{min}(\omega)\kappa_{min}(\omega_p)}.\label{split4}
		\end{eqnarray}
		As mentioned before, we split the concerned error into four parts as in \ref{split error}. From (\ref{split1}) and (\ref{split4}), one have
		\begin{eqnarray*}
			\|e_1\|_{L^2}^2\leq \dfrac{C\epsilon^2}{\kappa_{min}(\omega)\kappa_{min}(\omega_p)}, \text{ and } \|e_4\|_{L^2}^2\leq \dfrac{C\epsilon^2}{\kappa_{min}(\omega)\kappa_{min}(\omega_p)}.
		\end{eqnarray*}
		
		Combining Theorem \ref{e_2} and \ref{e_3}, we have (\ref{final_1}).
	\end{proof}

	\section{Numerical experiments}
	\begin{figure}[!htbp]
		\subfigure[$log(\kappa_1)$.]
		{ \includegraphics[width=0.45\textwidth]{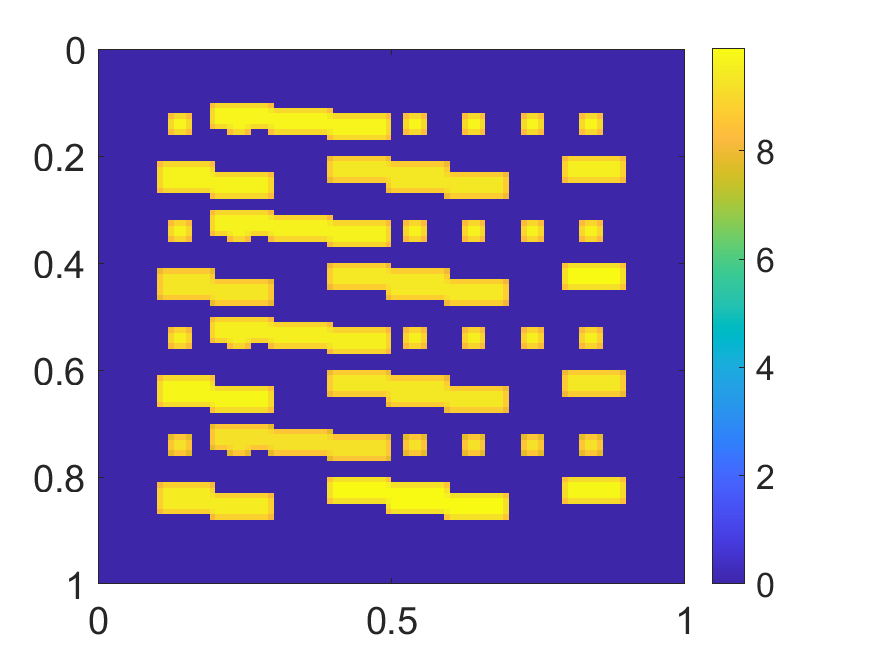}\label{fig:model1}}
		\subfigure[$log(\kappa_2)$.]
		{ \includegraphics[width=0.45\textwidth]{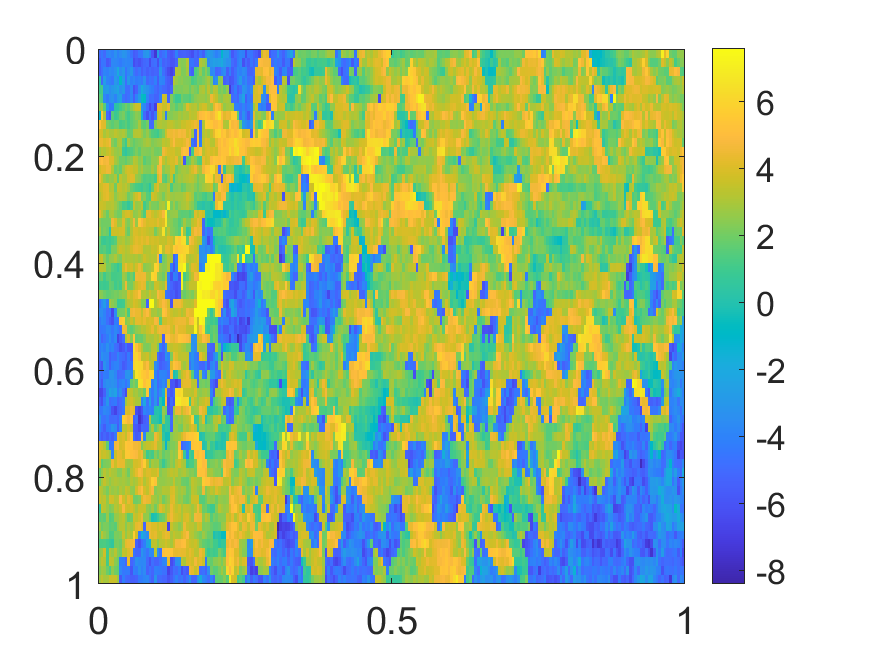}\label{fig:model2}}
		\caption{Display of fields in natural logarithm; Left: high-contrast model ; Right: SPE model}.\label{fig:model}
	\end{figure}
	In this section, we show the performance of our method. As shown in subsection \ref{sect:GMsFEM-POD}, we need to apply GMsFEM on mean permeability field. We consider two typical highly heterogeneous permeability fields to serve as mean permeability fields, whose natural logarithms are displayed in Figure \ref{fig:model}. Note that $\kappa_1$ is a deterministic permeability field with a high contrast of $10^4$ and resolution of $100\times 100$. $\kappa_2$ (see Figure \ref{fig:model2}) is the last layer of the commonly used SPE comparative solution project (SPE10), which is effective in assessing upscaling and multiscale methods. In spite of the fact that no distinct separate channels can be detected, it is highly heterogeneous and of high contrast with the magnitude approximately $10^7$. Besides, it has resolution $60\times 220$. To evaluate the accuracy  of our proposed method, we consider two types of relative errors defined below, i.e. energy error $e_a(t)$ and $L^2$ error $e_{L^2}(t)$ at time $t$. The reference solutions $u_{h}(\mathbf{x},t)$ is obtained by FEM with fine mesh $100\times 100$ for $\kappa_1$ and $60\times 220$ for $\kappa_2$.  $u_{appro}$ is the corresponding approximation solution. For coarse mesh, we use $10\times 10$ for $\kappa_1$ and $6\times 22$ for $\kappa_2$. Moreover, the final time $T$ is set to be 1.
	\begin{eqnarray}
		e_a(t)=\dfrac{\sqrt{\int_{\Omega}\kappa(\mathbf{x},\omega)|\nabla u_{appro}(\mathbf{x},t)-\nabla u_{h}(\mathbf{x},t)|^2 }}{\sqrt{\int_{\Omega}\kappa(\mathbf{x},\omega)|\nabla u_{h}(\mathbf{x},t)|^2 }}.\\
		e_{L^2}(t)=\dfrac{\sqrt{\int_{\Omega}|u_{appro}(\mathbf{x},t)-u_{h}(\mathbf{x},t)|^2}}{\sqrt{\int_{\Omega}|u_{h}(\mathbf{x},t)|^2 }}.
		\label{err_num}
	\end{eqnarray}
	\begin{table}[htbp!]
		\centering
		\begin{tabular}{c|c|c}
			Method     &  Dimension of system & Solving time(sec)\\
			\hline \hline
			FEM  &    13481             &       0.0687     \\ 
			\hline
			GMsFEM(2+3)  &   805        &       0.01383    \\
			\hline
			GMsFEM(3+3)  &  966         &       0.02317     \\
			\hline
			GMsFEM(4+3)  &  1127        &       0.03242     \\
			\hline
			GMsFEM(8+0)  &  1288        &       0.04301     \\
			\hline
			GMsFEM-POD(5)       &   5          &       0.00031      \\
			\hline
			GMsFEM-POD(10)      &   10         &       0.00054      \\
			\hline
			GMsFEM-POD(15)      &   15         &       0.0001      \\
			\hline
			GMsFEM-POD(20)      &   20         &       0.0011      \\
			\hline
			GMsFEM-POD(25)      &   25         &       0.00012      \\
			\hline
		\end{tabular}
		\caption{Computation time(sec): per single solve of the resulting linear system.}
		\label{tab:solving time}
	\end{table}
	We first show the computational time per single solve of the resulting linear system with different methods in Table \ref{tab:solving time}. In particular, we show the computation time of solving (\ref{fine sol stoc}), (\ref{ms sol stoc}), and (\ref{snap sol stoc}) with FEM, GMsFEM and GMsFEM-POD, respectively. For GMsFEM,
	there are two separate steps in construction of offline space represented in Table \ref{online algo1} and \ref{online algo2}. Here we use $``A+B"$ to denote the case that $A$ basis functions are used in each local neighborhood in offline stage 1 while $B$ additional basis functions are incorporated in stage 2. From Table \ref{tab:solving time}, we can see the computation time in GMsFEM-POD system is much shorter than GMsFEM and FEM system, which clearly shows the efficiency of the proposed method.
	
	Next we consider the effects of important components of the proposed method in the following part.
	\subsection{Influence of the parameters in KLE.}
\begin{figure}[!htbp]
		\centering
		\subfigure[High-contrast model.]
		{ \includegraphics[width=0.4\textwidth]{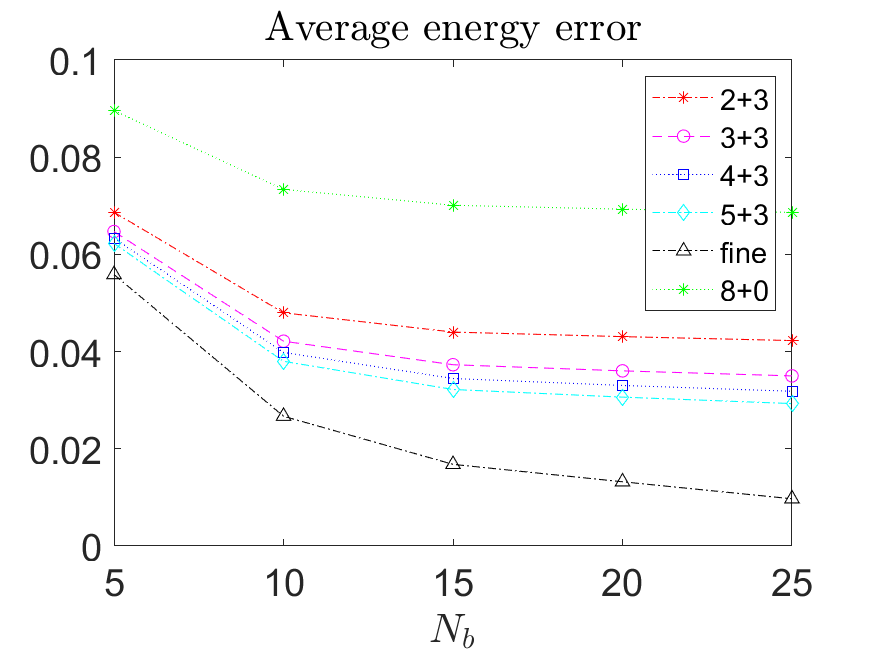}\label{fig:sigma2_2_0.5_0.5_a}}
		\subfigure[SPE model.]
		{ \includegraphics[width=0.4\textwidth]{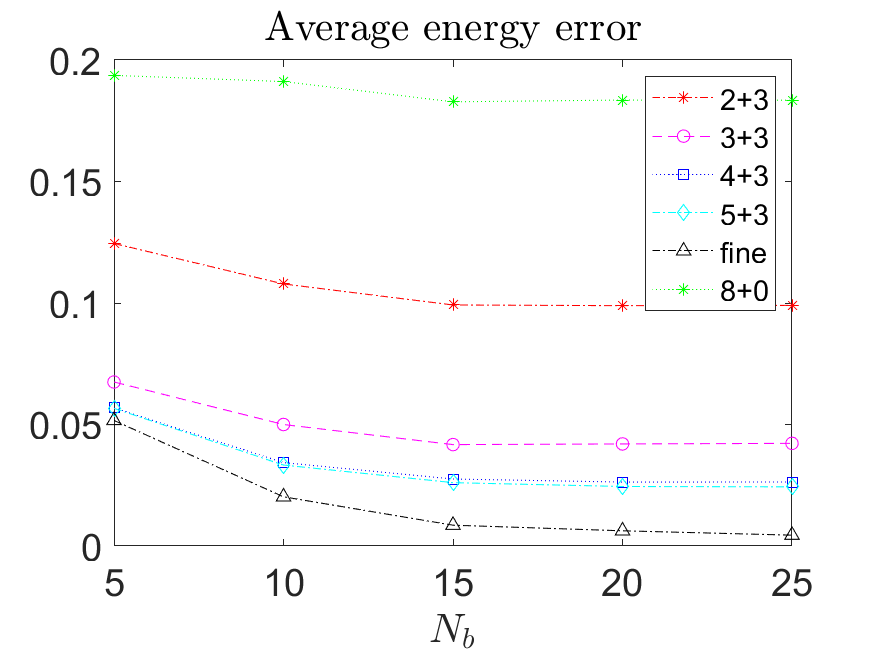}\label{fig:sigma2_2_0.5_0.5_b}}
		\subfigure[High-contrast model.]
		{ \includegraphics[width=0.4\textwidth]{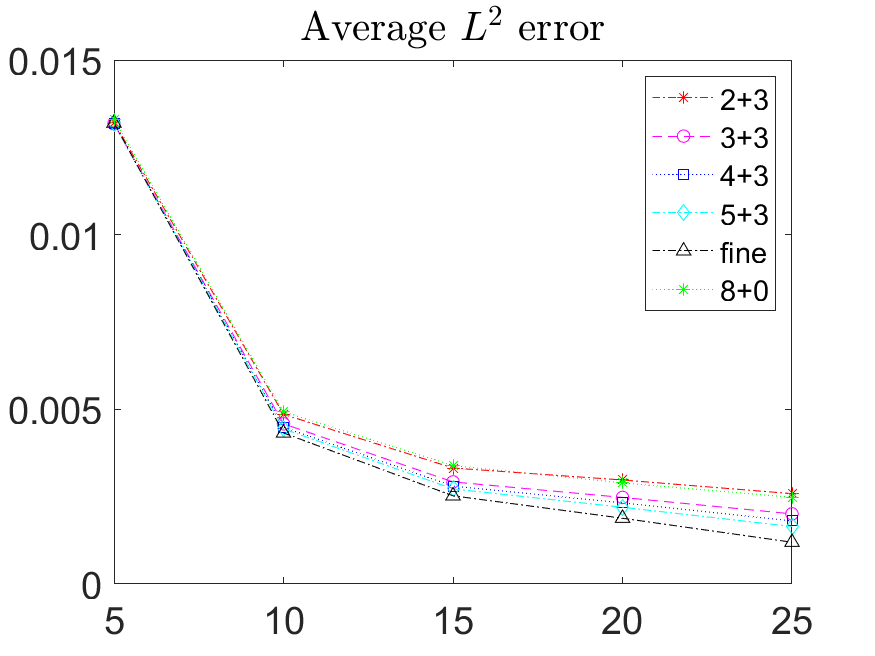}}
		\subfigure[SPE model.]
		{ \includegraphics[width=0.4\textwidth]{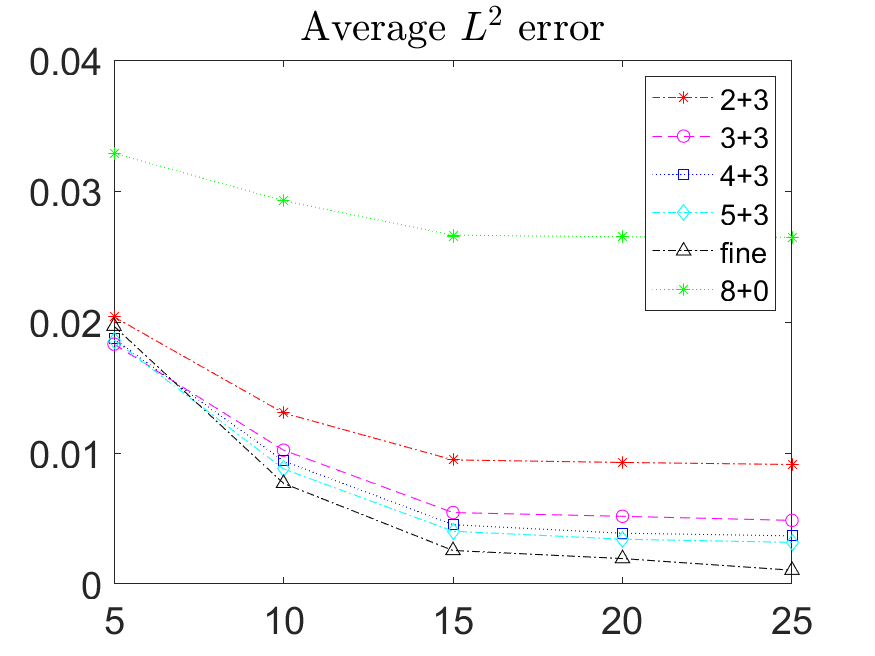}}
		\subfigure[High-contrast model.]
		{ \includegraphics[width=0.4\textwidth]{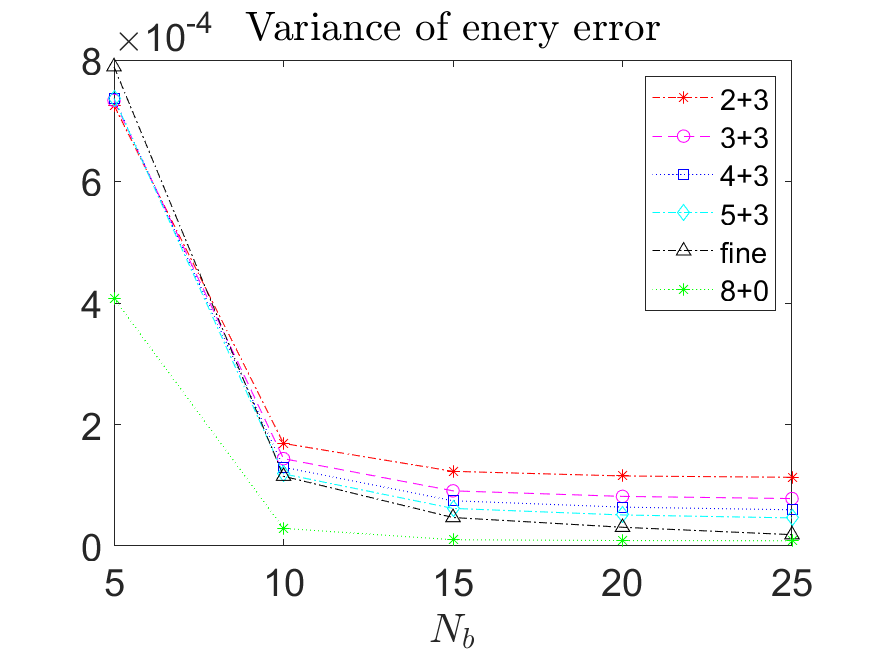}}
		\subfigure[SPE model.]
		{ \includegraphics[width=0.4\textwidth]{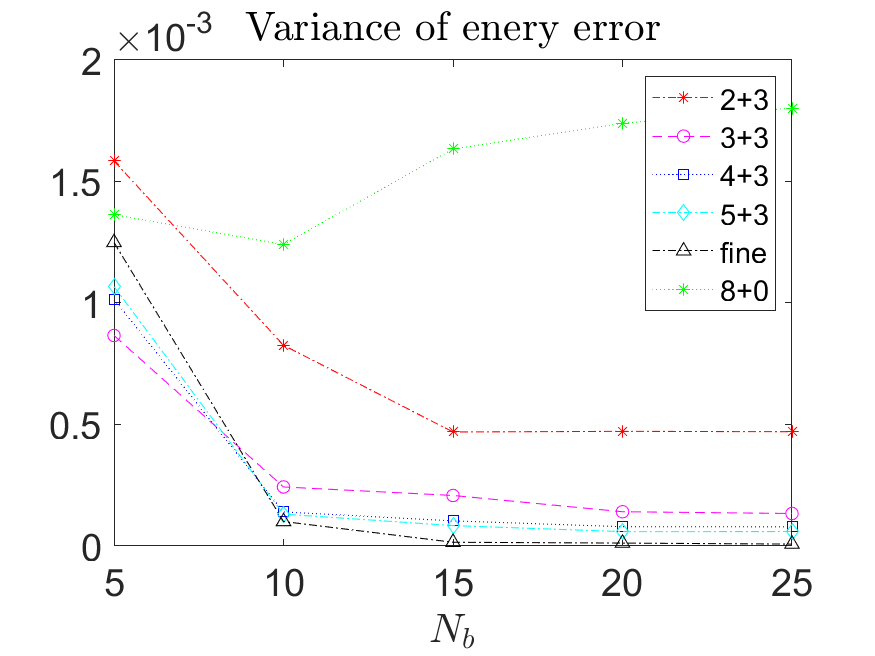}}
		\subfigure[High-contrast model.]
		{ \includegraphics[width=0.4\textwidth]{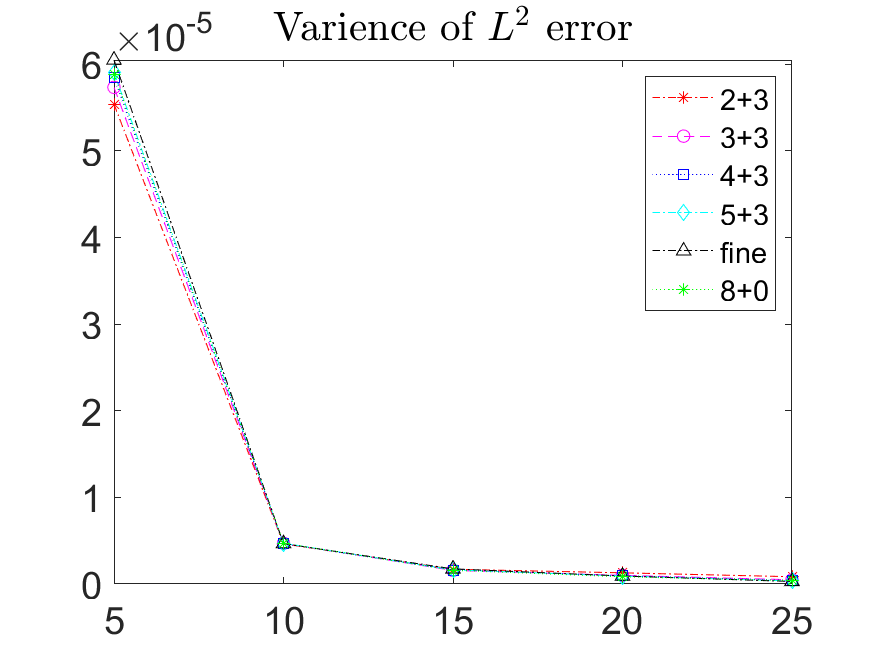}}
		\subfigure[SPE model.]
		{ \includegraphics[width=0.4\textwidth]{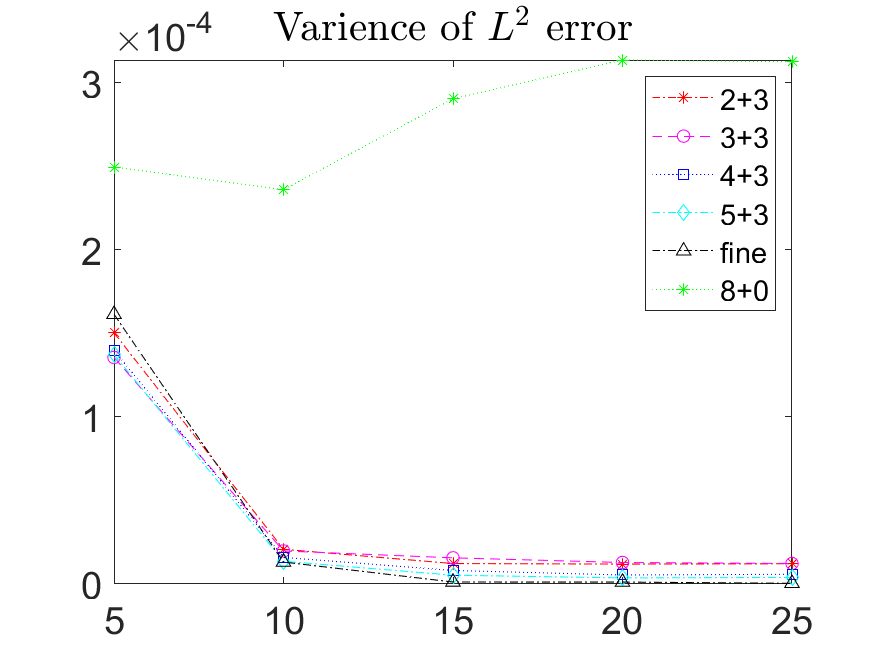}}
		\caption{Mean and variance of energy and $L^2$  errors (of 100 samples) at final time $T$ under $\sigma^2=2, \eta_1=0.5, \eta_2=0.5$.}\label{fig:sigma2_2_0.5_0.5}
	\end{figure}

\begin{figure}[!htbp]
		\centering
		\subfigure[High-contrast model.]
		{ \includegraphics[width=0.4\textwidth]{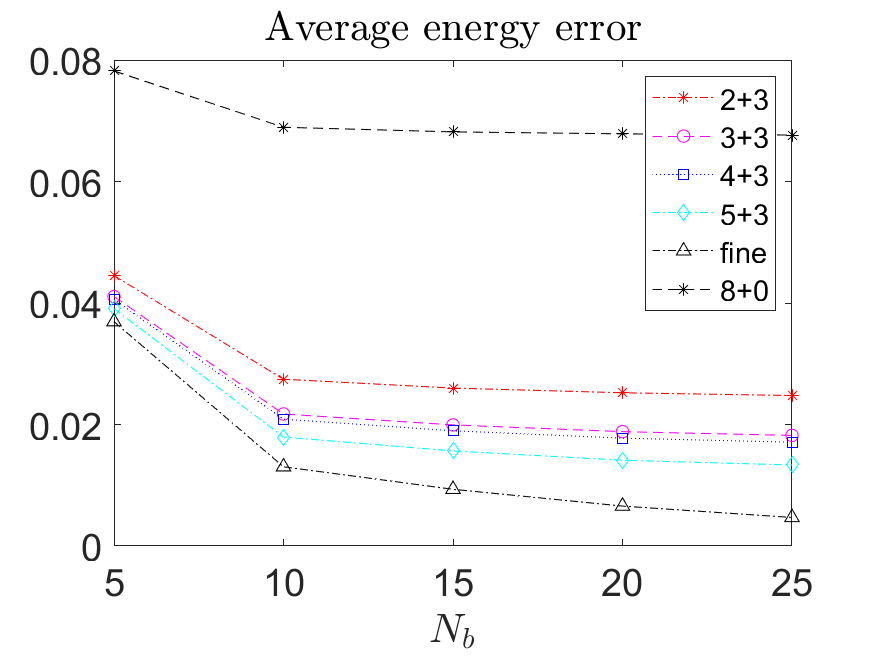}\label{fig:sigma2_1_0.5_0.5_a}}
		\subfigure[SPE model.]
		{ \includegraphics[width=0.4\textwidth]{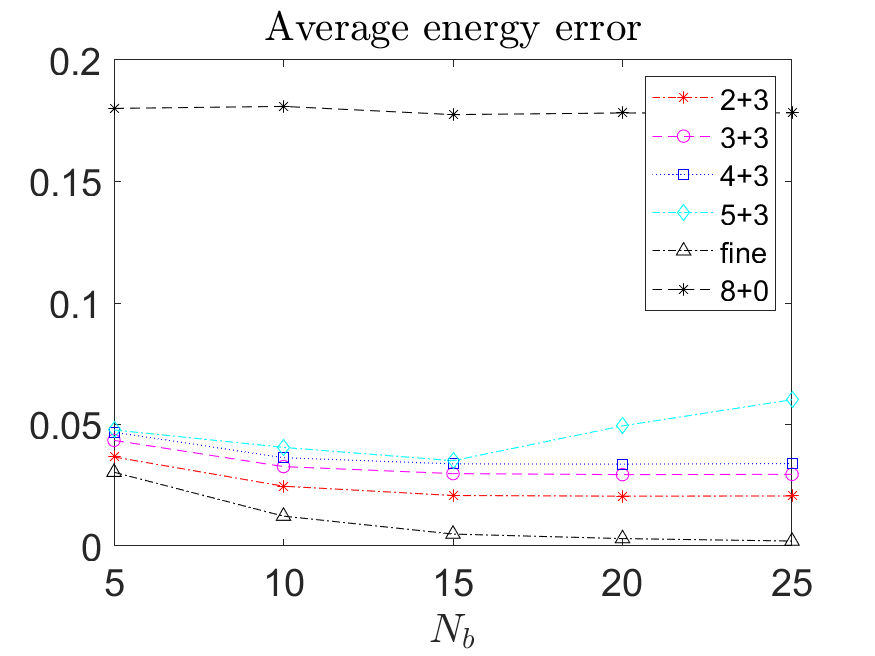}\label{fig:sigma2_1_0.5_0.5_b}}
		\subfigure[High-contrast model.]
		{ \includegraphics[width=0.4\textwidth]{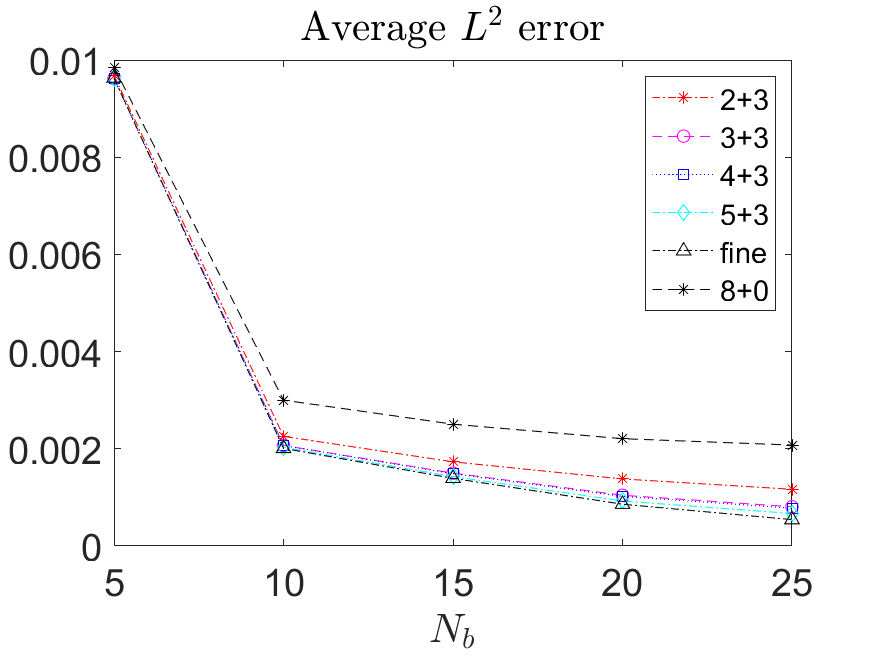}}
		\subfigure[SPE model.]
		{ \includegraphics[width=0.4\textwidth]{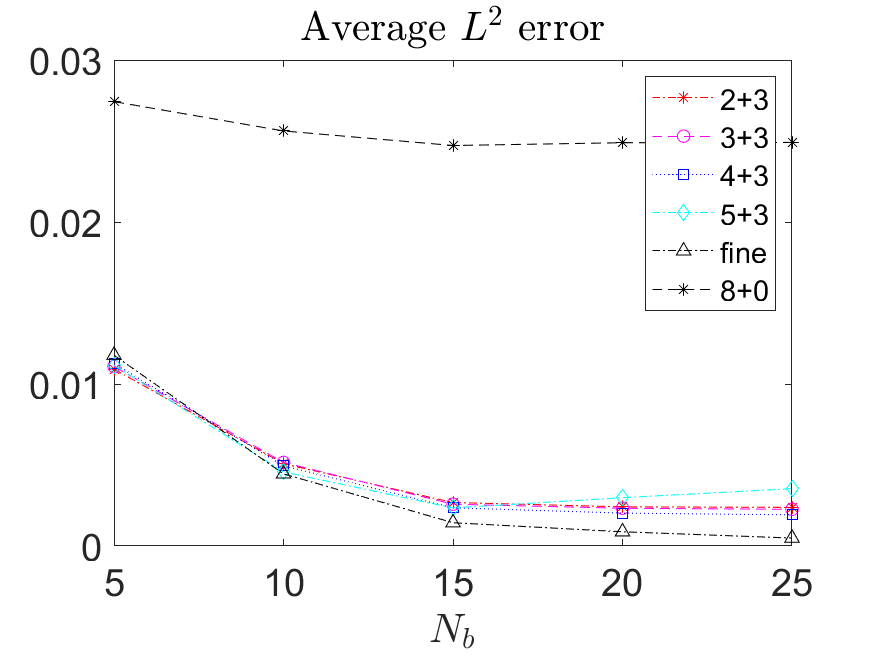}}
		\subfigure[High-contrast model.]
		{ \includegraphics[width=0.4\textwidth]{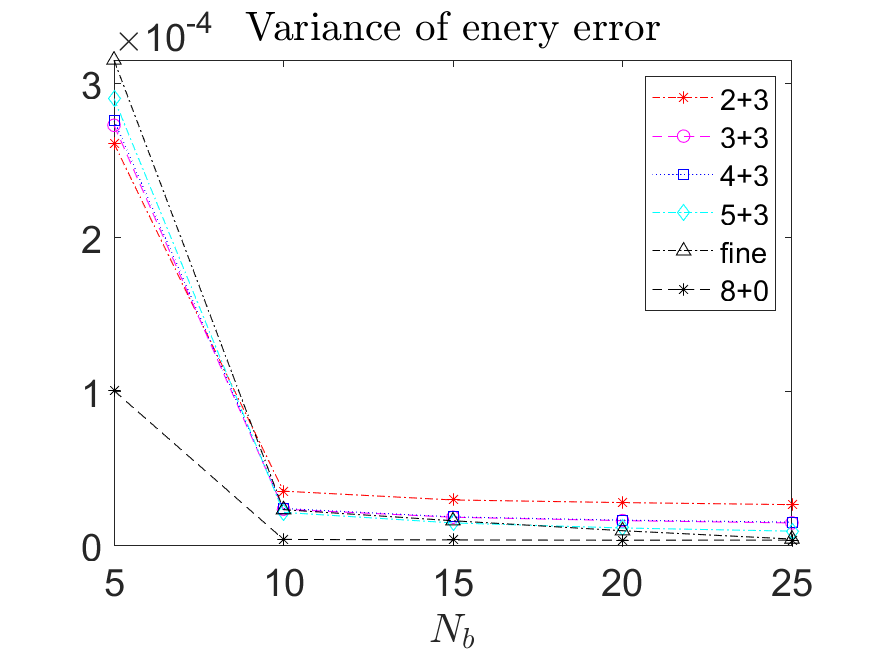}}
		\subfigure[SPE model.]
		{ \includegraphics[width=0.4\textwidth]{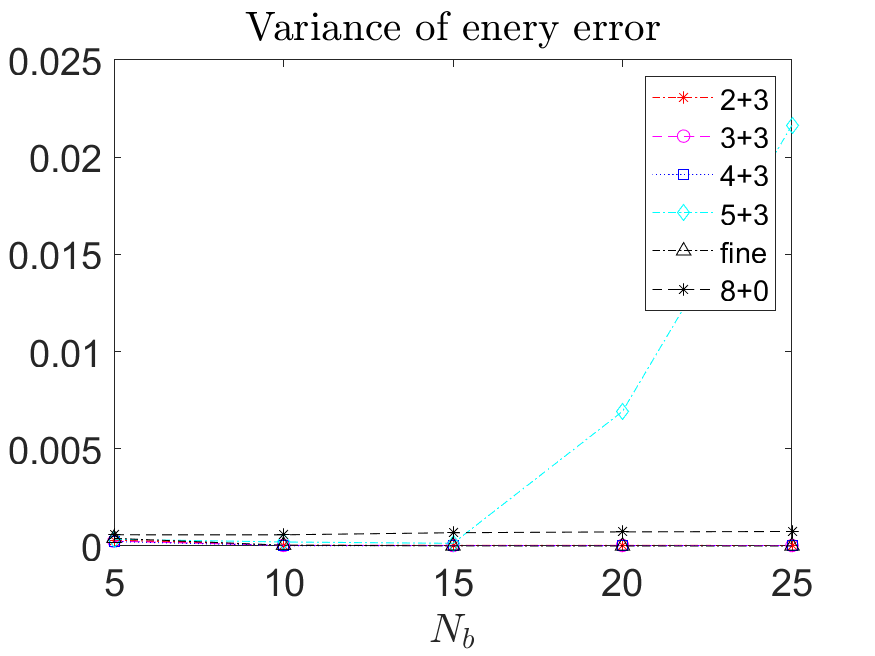}}
		\subfigure[High-contrast model.]
		{ \includegraphics[width=0.4\textwidth]{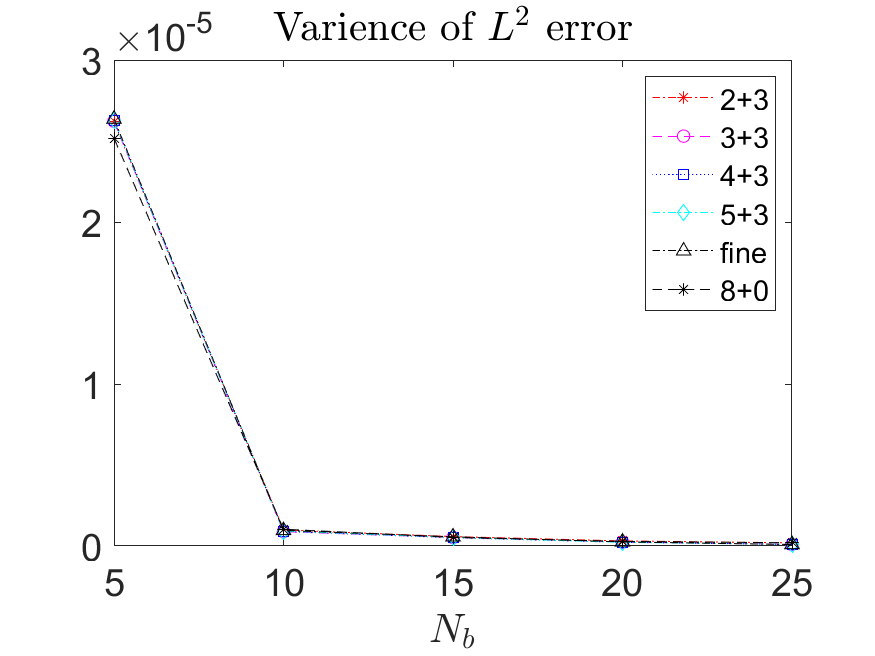}}
		\subfigure[SPE model.]
		{ \includegraphics[width=0.4\textwidth]{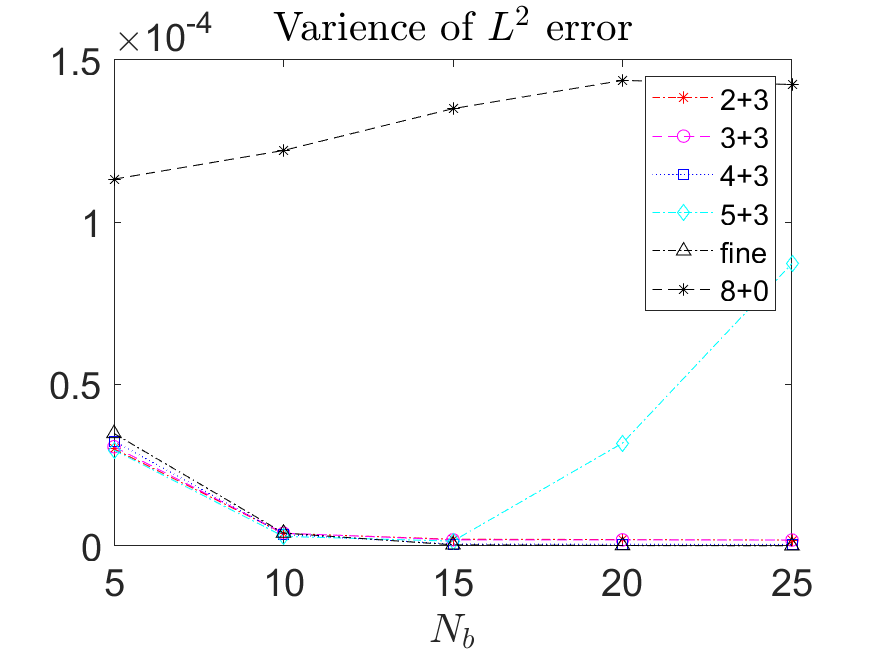}}
		\caption{Mean and variance of energy and $L^2$  errors (of 100 samples) at final time $T$ under $\sigma^2=1, \eta_1=0.5, \eta_2=0.5$.}\label{fig:sigma2_1_0.5_0.5}
	\end{figure}

\begin{figure}[!htbp]
		\centering
		\subfigure[High-contrast model.]
		{ \includegraphics[width=0.4\textwidth]{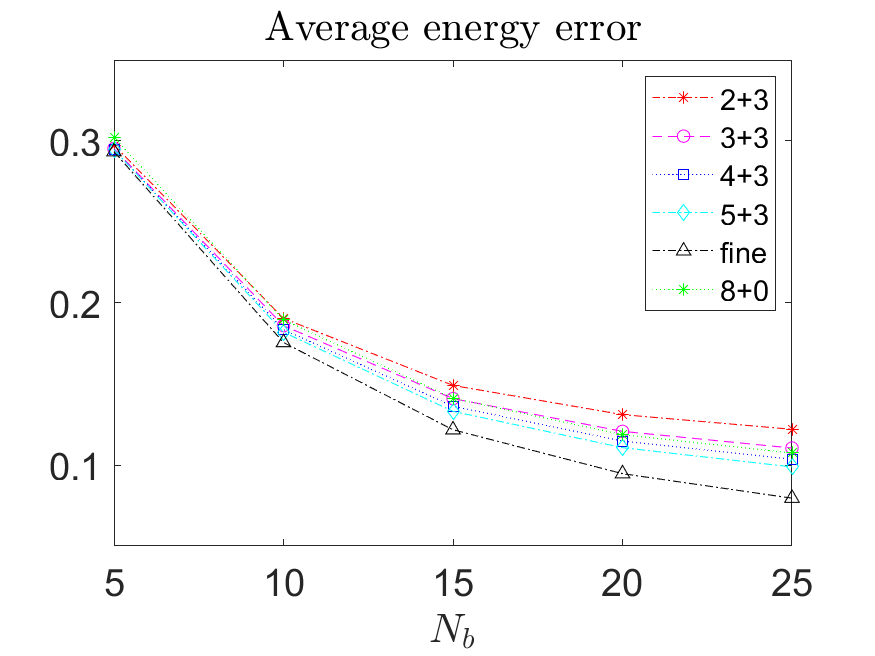}}
		\subfigure[SPE model.]
		{ \includegraphics[width=0.4\textwidth]{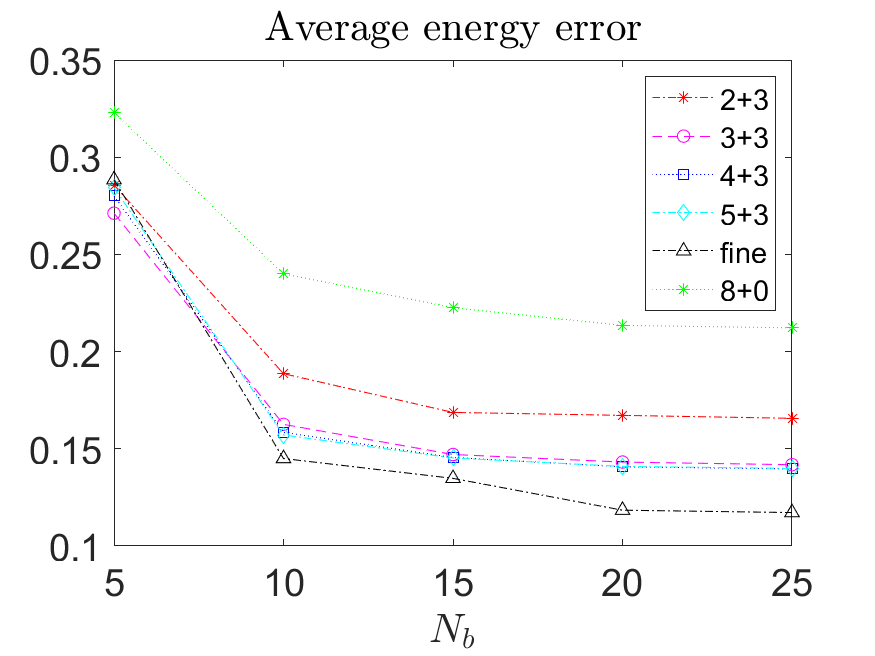}}
		\subfigure[High-contrast model.]
		{ \includegraphics[width=0.4\textwidth]{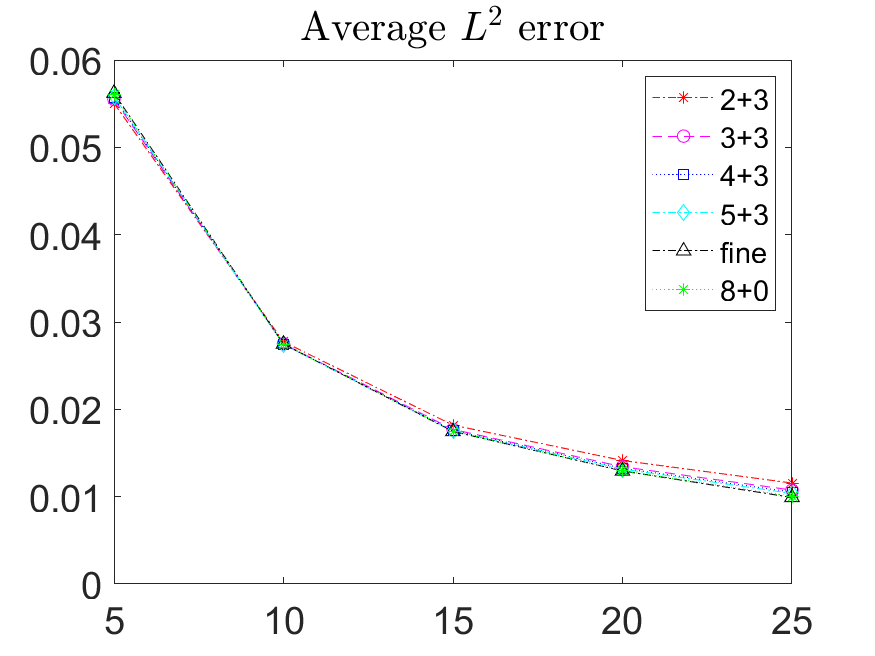}}
		\subfigure[SPE model.]
		{ \includegraphics[width=0.4\textwidth]{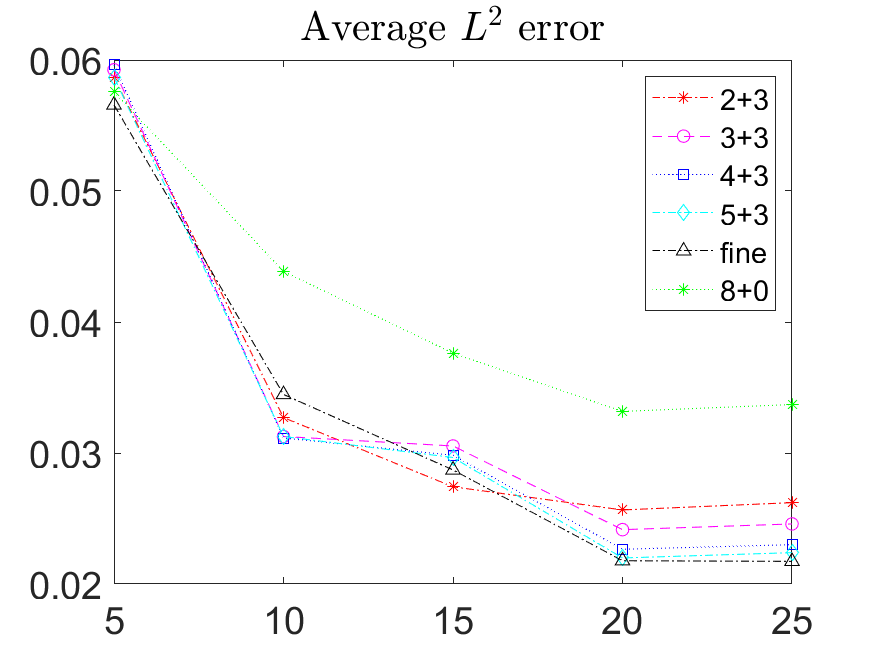}}
		\subfigure[High-contrast model.]
		{ \includegraphics[width=0.4\textwidth]{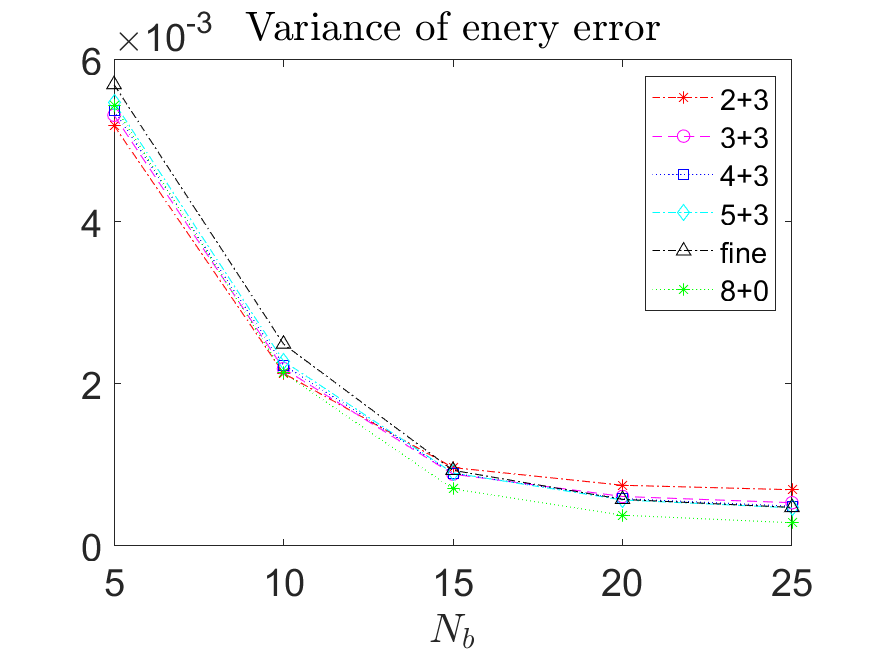}}
		\subfigure[SPE model.]
		{ \includegraphics[width=0.4\textwidth]{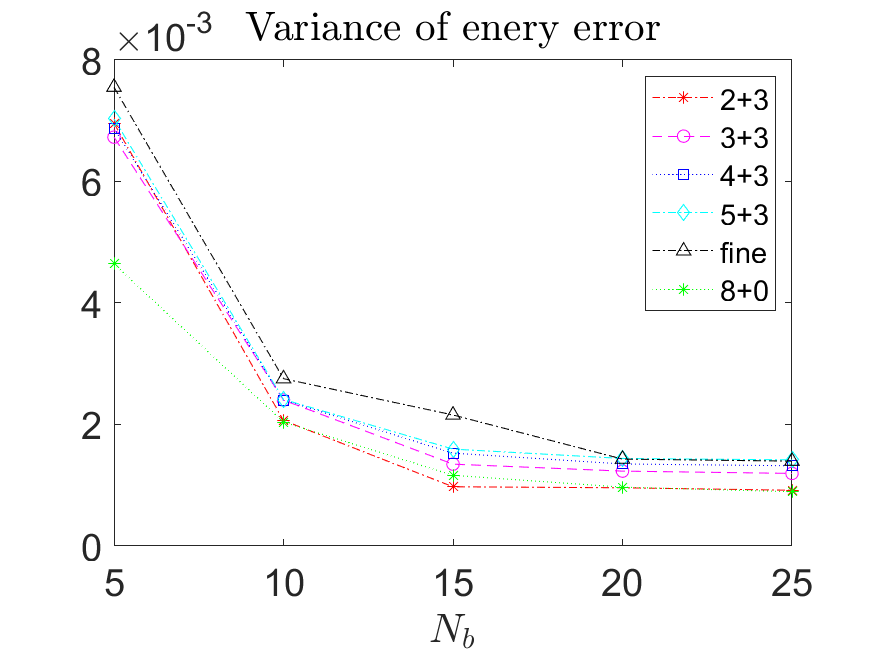}}
		\subfigure[High-contrast model.]
		{ \includegraphics[width=0.4\textwidth]{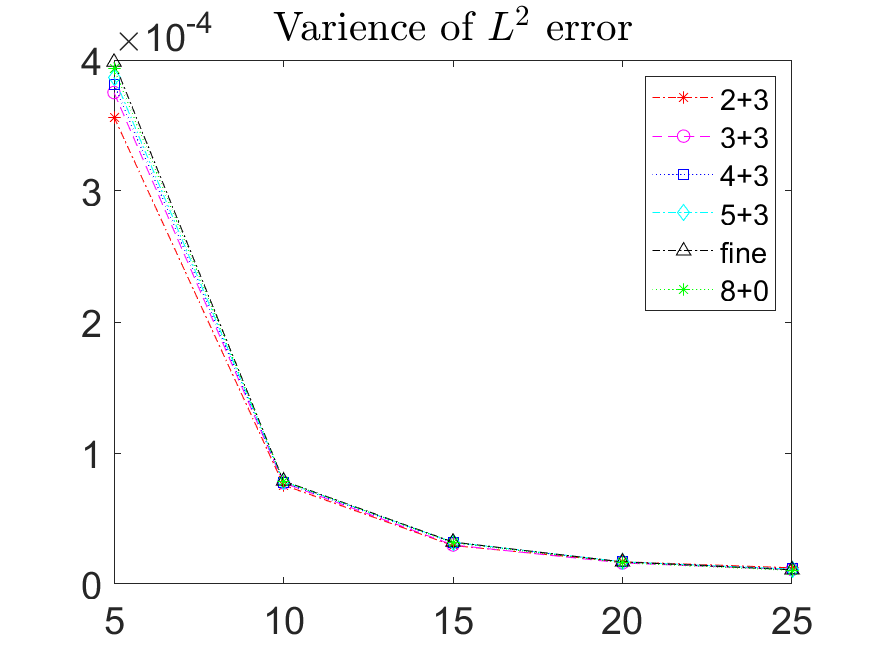}}
		\subfigure[SPE model.]
		{ \includegraphics[width=0.4\textwidth]{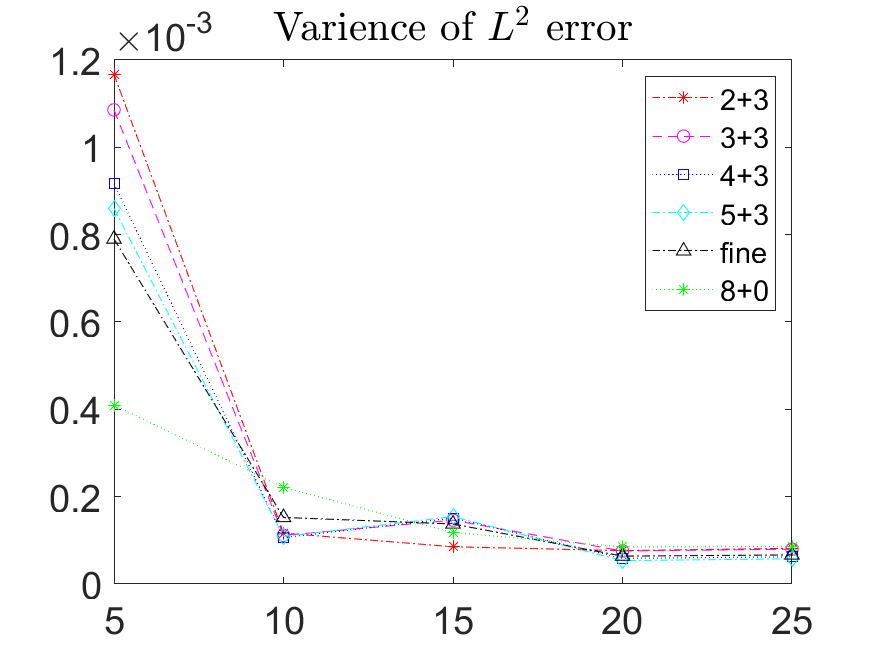}}
		\caption{Mean and variance of energy and $L^2$  errors (of 100 samples) at final time $T$ under $\sigma^2=1, \eta_1=0.5, \eta_2=0.1$.}\label{fig:sigma2_1_0.5_0.1}
	\end{figure}

	In this subsection, we mainly study the effects of the multiscale space dimension in offline stage 1 and 2 under different choices of covariance matrices and two types of mean permeability fields (in Figure \ref {fig:model}). In particular, we consider three cases of $\sigma^2$, $\eta_1$ and $\eta_2$ as follows.
	\begin{enumerate}
		\item Figure \ref{fig:sigma2_2_0.5_0.5}: $\sigma^2=2$, $\eta_1=0.5$, $\eta_2=0.5$.
		\item Figure \ref{fig:sigma2_1_0.5_0.5}: $\sigma^2=1$, $\eta_1=0.5$, $\eta_2=0.5$.
		\item Figure \ref{fig:sigma2_1_0.5_0.1}: $\sigma^2=1$, $\eta_1=0.5$, $\eta_2=0.1$.
	\end{enumerate}
	We first study the influence of the number of global POD basis and local multiscale 
	basis.
	In Figure \ref{fig:sigma2_2_0.5_0.5} - \ref{fig:sigma2_1_0.5_0.1}, we display the dynamics of interested quantities (mean and variance of 100 samples) of the energy errors and  $L_2$ errors associated with high-contrast model and SPE model.  In each subfigure, it displays errors with different numbers of POD basis and multiscale basis. 
	To investigate the influence of using less accurate global snapshots (multiscale solutions), we also 
	show the errors of using fine-grid solutions as the global snapshots for POD basis construction in each 
	test case for comparison, where we use ``fine" to denote this case. Indeed we can treat this ``fine"  error as the limiting error of our local-global method.
	From these error plots we first observe that with the same number of multiscale basis functions, both the average and variance of the energy errors and  $L_2$ errors decay as more POD basis are used. However, both types of errors are tending to be stable as the number of POD basis exceeds 15. Second observation is that with the number of POD basis fixed, the errors decrease if more local multiscale basis are used in the training stage. Another important observation is that the effect of residual-driven basis is impressive as one compares the $``2+3"$ and $``8+0"$ case. In particular, in Figure \ref{fig:sigma2_2_0.5_0.5_a}, the average energy error ends with about $7\%$ in the case $``8+0"$ while for the case $``2+3"$, corresponding error drops below $4\%$. In a word, less local basis functions can be used if one includes residual-driven basis to obtain desired accuracy.
	Last observation is that using the globally defined multiscale snapshots does not bring too much additional error in the POD method by comparing the error of the cases like ``5+3'' and ``fine'' case, where relatively small difference is observed between these two cases.

\begin{figure}[!htbp]
		\centering
		\subfigure[$2+3$ basis.]
		{ \includegraphics[width=0.45\textwidth]{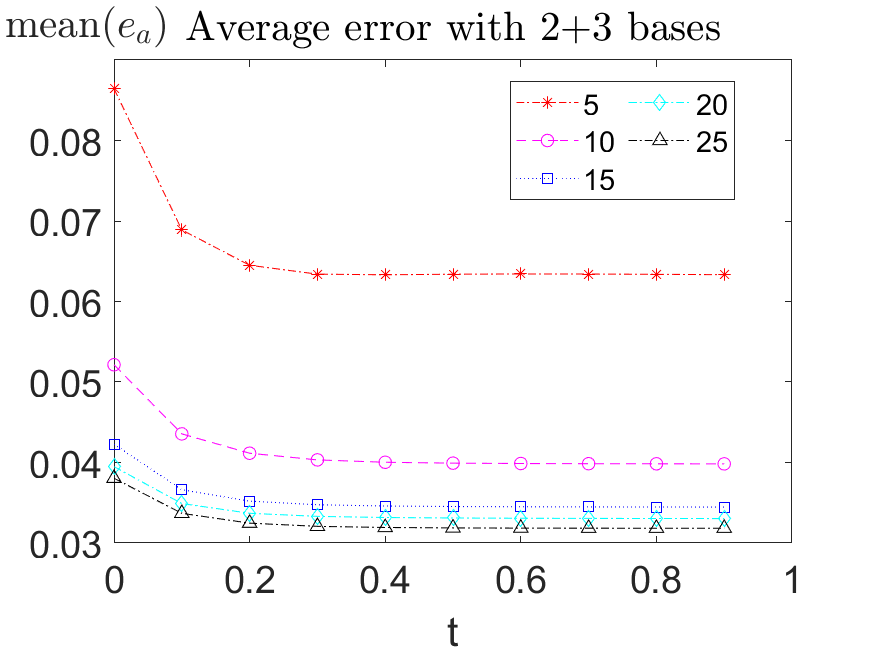}\label{fig:sigma2_1_0.5_0.1_dynamics_a}}
		\subfigure[$3+3$ basis.]
		{ \includegraphics[width=0.45\textwidth]{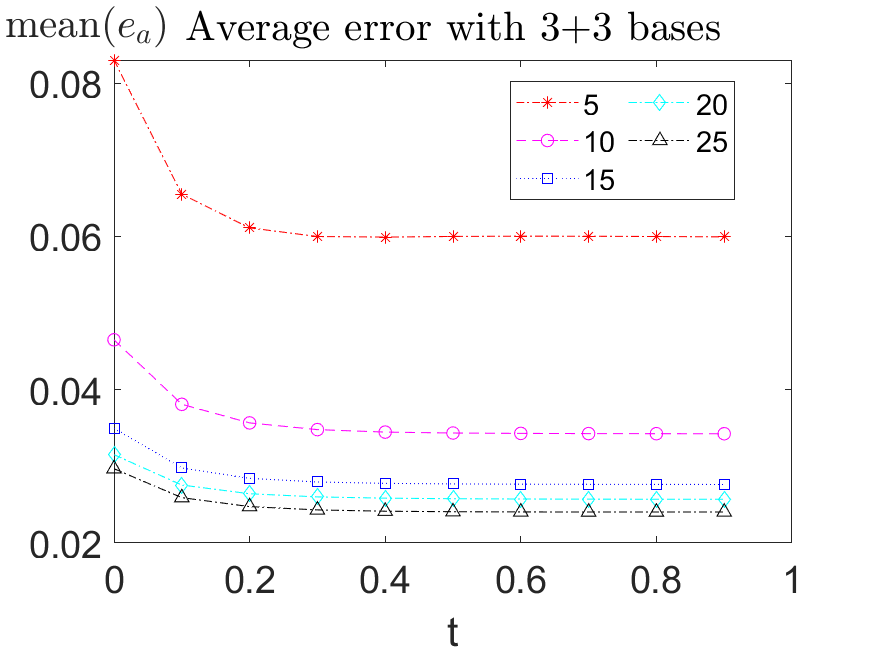}}
		\subfigure[$4+3$ basis.]
		{ \includegraphics[width=0.45\textwidth]{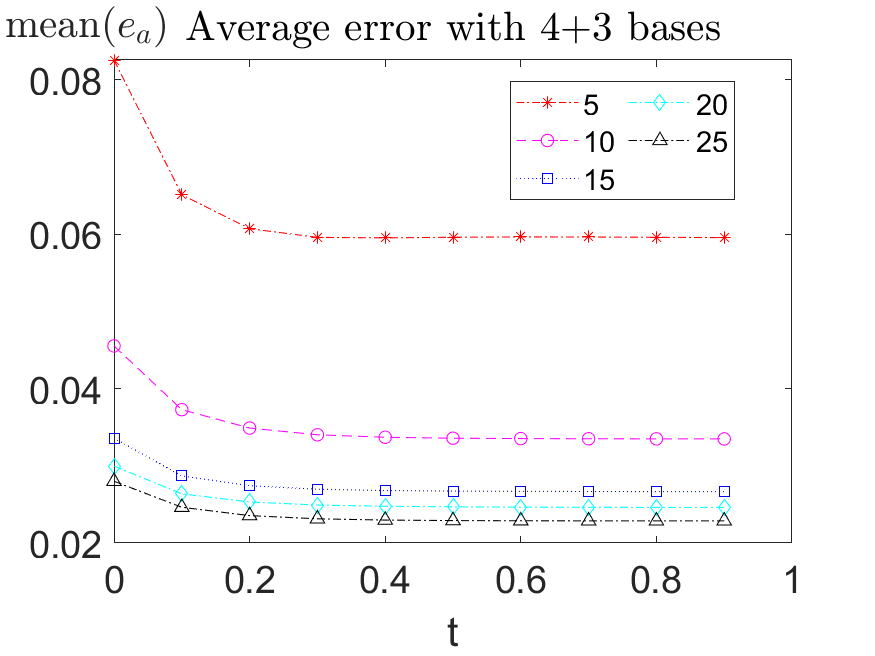}}
		\subfigure[$8+0$ basis.]
		{ \includegraphics[width=0.45\textwidth]{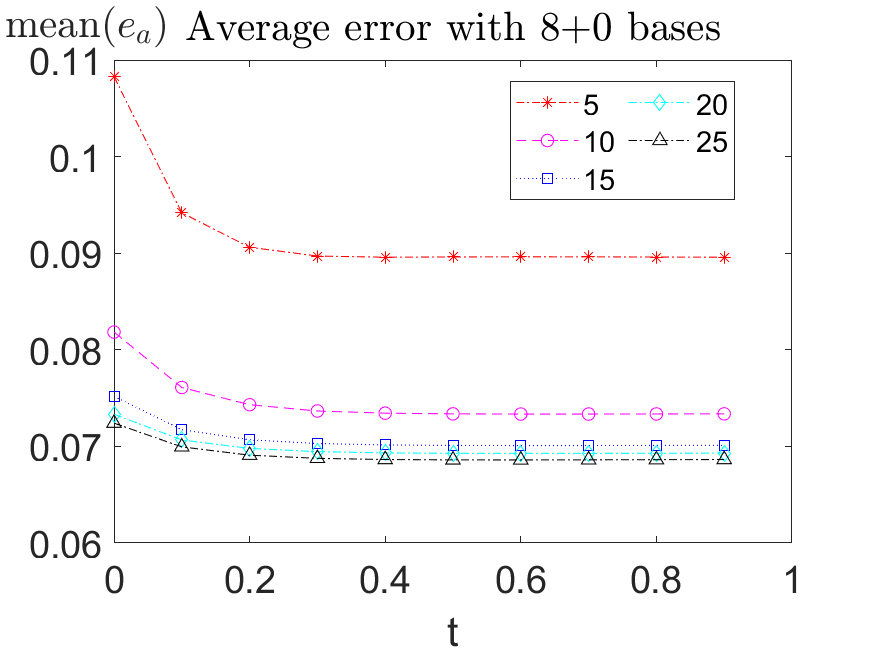}}
		\caption{Dynamics of average energy errors (of 1000 samples) with different numbers of basis under $\sigma^2=2, \eta_1=0.5, \eta_2=0.5$.}
		\label{fig:sigma2_1_0.5_0.1_dynamics}
	\end{figure}

	We then investigate the influence of using different numbers of POD basis in detail. Here, we fix $\sigma^2=2, \eta_1=0.5, \eta_2=0.5$. And we use high-contrast model. Rather than simply displaying the error of final states as in Figure \ref{fig:sigma2_2_0.5_0.5}, we  show the dynamics of average energy error (of 100 samples) with respect to time in Figure \ref{fig:sigma2_1_0.5_0.1_dynamics}. In spite of using various numbers of POD basis functions, the dynamics of the errors behave in a similar manner. For example, in Figure \ref{fig:sigma2_1_0.5_0.1_dynamics_a}, we compare using 5, 10, 15, 20, 25 POD basis functions while the number of offline local basis functions is chosen to be the same ($``2+3"$). We see that in general, the error decreases as the time advances before the first half of total time period and keeps steady in the rest time period. In Figure \ref{fig:sigma2_1_0.5_0.1_dynamics_a}, there is a sharp decline of error  from using 5 POD basis functions to the case where 15 basis are used, in particular, from around $8.5\%$ to $3.5\%$. However, the improvement is less significant as further enrichment is performed in the POD basis functions. For the other three cases, i.e. $``3+3"$, $``4+3"$, $``8+0"$, similar changing behaviors are displayed. Therefore, one could see the improvement of adding POD basis is evident when the initial number is small. Furthermore, a few numbers of POD basis (10 in our example) are sufficient to get relatively high accuracy. However, similar to the case of enriching local multiscale basis, once the number of basis exceeds certain bounds, the corresponding error reduction is not profound.

	\subsection{Comparison of two offline enrichment algorithms.}
	In this subsection, we compare two offline enrichment algorithms and simply use the SPE10 model as test model in step 1 of  both algorithms. We present the average energy errors (of 100 samples) in three steps with different numbers of multiscale basis in Table \ref{tab:algo1_sigma2_5} to \ref{tab:algo2_sigma2_1}. Here we use  $``A+0+B+C+D"$ to denote the case where $A$ basis are used in offline stage 1 with mean permeability field while $B$, $C$, $D$ residual-driven basis functions associated with different chosen samples  incorporated hierarchically in stage 2. For all the three steps, the reference solutions are solved with standard finite element methods on fine mesh, corresponding to (\ref{fine sol}). In the first and second steps, we seek solutions to (\ref{ms sol}) and (\ref{ms sol stoc}) respectively with multiscale basis. In the third step, we seek solutions to (\ref{snap sol stoc}) with POD basis. In the first step, the error is corresponding to the mean permeability field. For the second and third steps, we compute mean energy error among 100 examples.  Moreover, we use 20 POD basis functions in third step.
	
	In Table \ref{tab:algo1_sigma2_5} and \ref{tab:algo1_sigma2_1}, we show the errors in three steps using algorithm 1 (Table \ref{tab:algo1}) for $\sigma^2=5$ and $\sigma^2=1$ respectively. We set $\eta_1=\eta_2=0.1$ for these two examples. The errors decay rapidly once residual-driven basis functions are used. In particular, as one can see from Table \ref{tab:algo1_sigma2_5}, errors decay from $30\%$ in $``5+0"$ to $1.4\%$ in $``2+3"$ in first step. However in the case of $\sigma^2=5$, the error decay in last two steps is not significant compared with the low-variance case. With $\sigma^2=1$, the error decreases sharply from over $30\%$ in $``5+0"$ to less than $10\%$ in $``2+3"$ in second and third steps. Moreover, we can see the accuracy in $``2+3"$ case is better than $``10+0"$ case, which shows the efficiency of residual-driven basis functions.
	
	In Table \ref{tab:algo2_sigma2_5} and \ref{tab:algo2_sigma2_1}, we shows the results using algorithm 2 (Table \ref{tab:algo2}). For step 2, we display errors in different offline enrichment levels. For example, in $``2+1+1+1"$ case, we split the enrichment process into three times. For each time, we add one basis function per local neighborhood and record the corresponding error after this enrichment. Comparing $``2+3"$ with $``2+1+1+1"$ case, the accuracy is further improved due to the splitting strategy in spite of same dimension of the equation system. In particular, in Table \ref{tab:algo2_sigma2_1}, the final error in $``2+1+1+1"$ case is about the half of the $``2+3"$ case and much smaller than $``5+0"$ case, which demonstrates the power of adopting splitting strategy in offline enrichment.
	
\begin{figure}[!htbp]
		\centering
		\subfigure[Reference solution in Step1.]
		{ \includegraphics[width=0.45\textwidth]{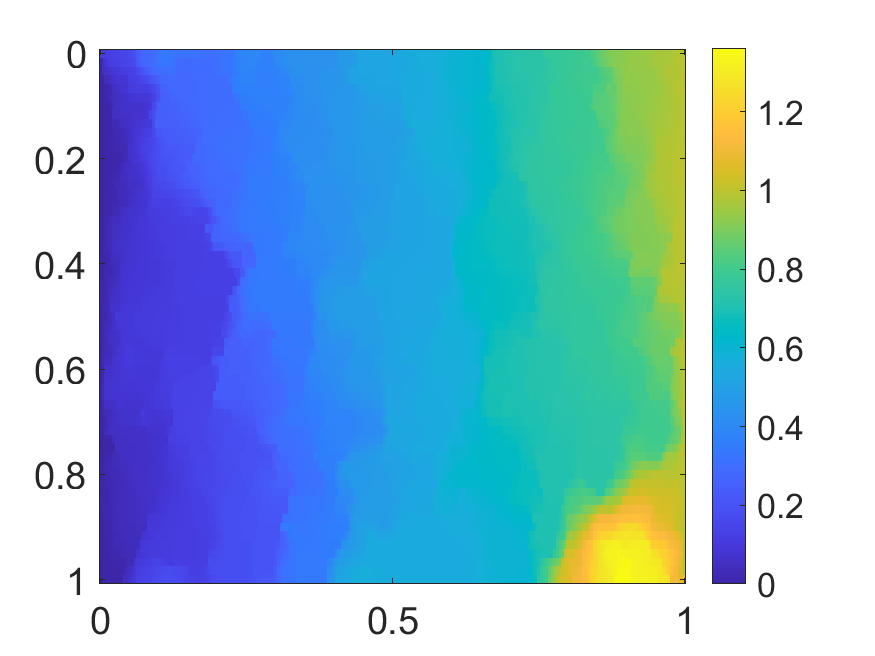}}
		\subfigure[Approximation in Step1.]
		{ \includegraphics[width=0.45\textwidth]{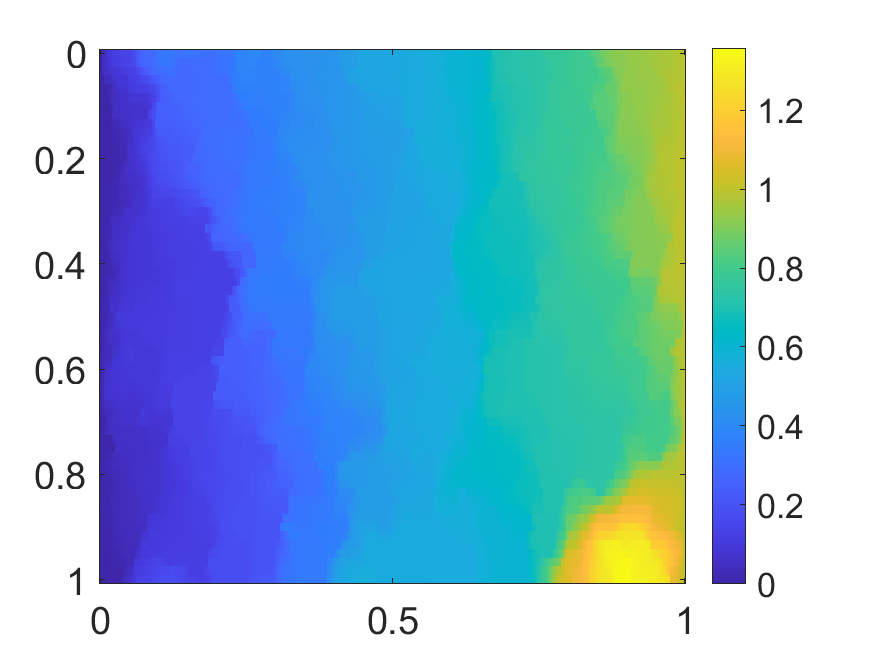}}
		\subfigure[Reference solution in Step2.]
		{ \includegraphics[width=0.45\textwidth]{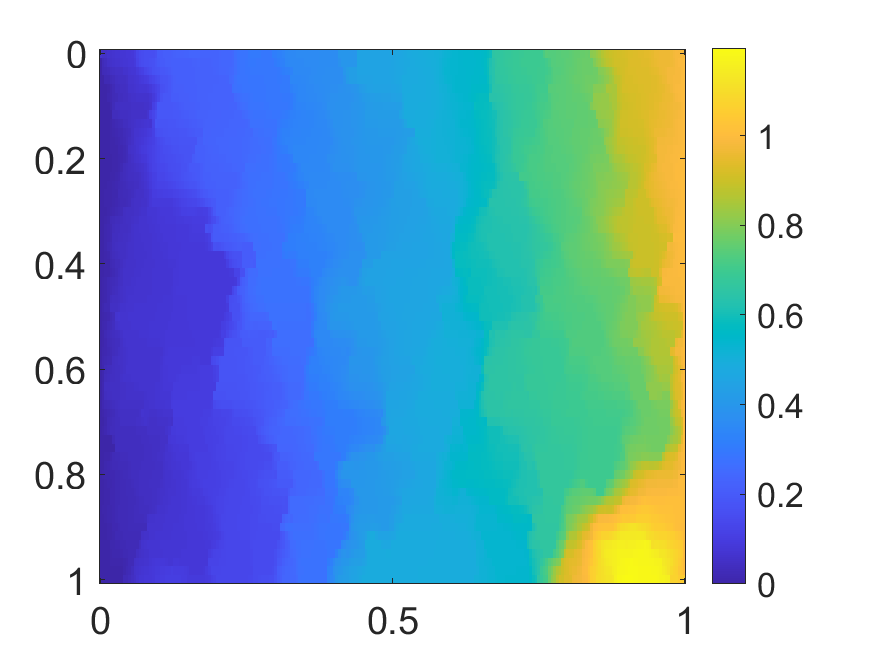}}
		\subfigure[Approximation in Step2.]
		{ \includegraphics[width=0.45\textwidth]{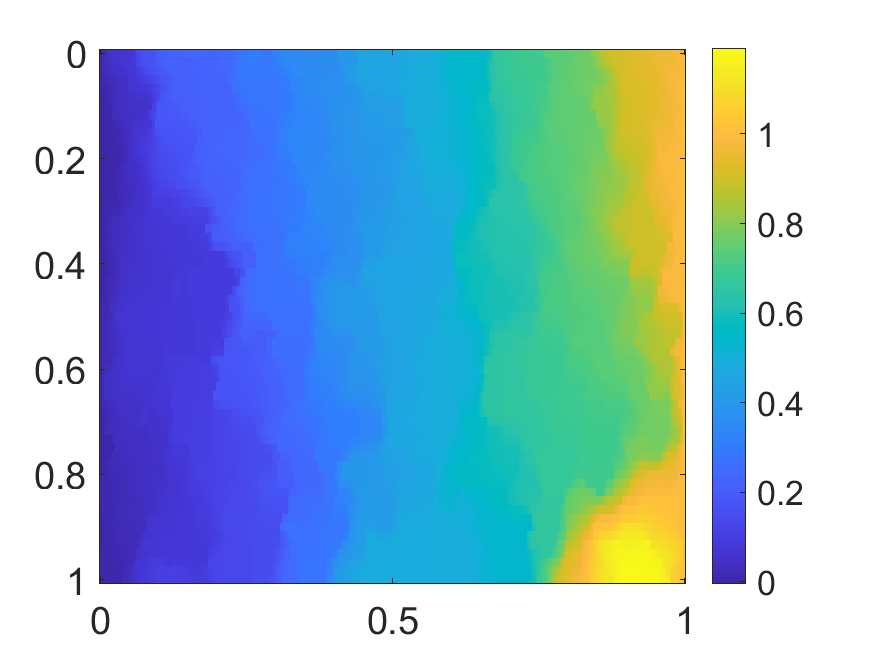}}
		\subfigure[Reference solution in Step3.]
		{ \includegraphics[width=0.45\textwidth]{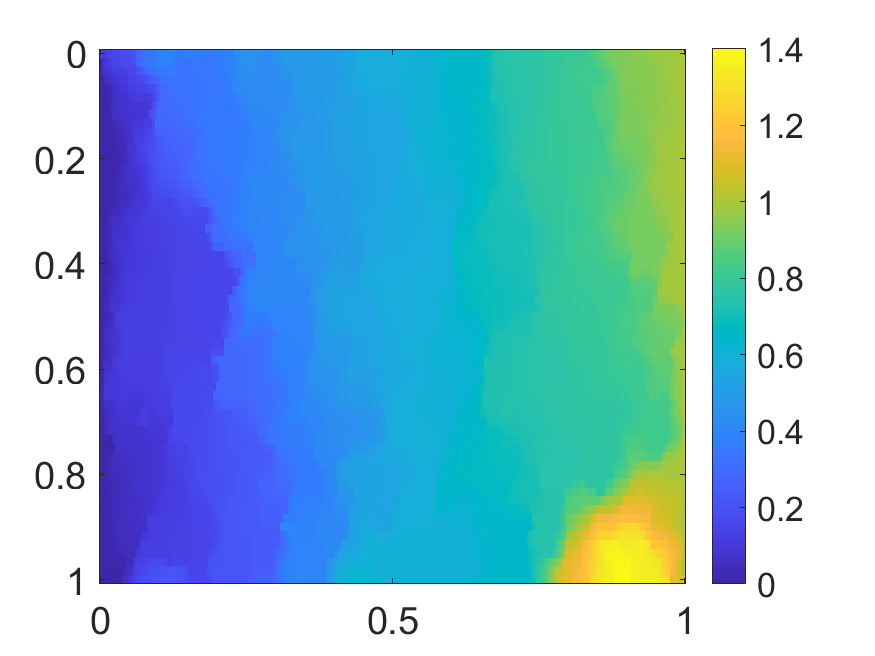}}
		\subfigure[Approximation in Step3.]
		{ \includegraphics[width=0.45\textwidth]{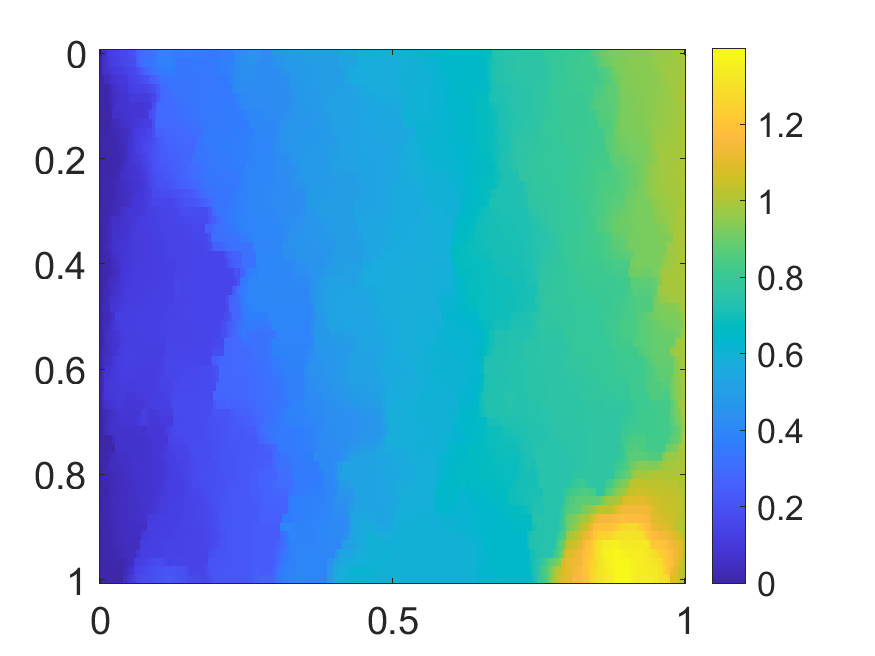}}
		\caption{Comparison of solutions of specific samples at time $t=T$ in three different steps in Algorithm 1. We use SPE model and set $\sigma^2=2$ and $\eta_1=\eta_2=0.5$.}\label{fig:algo1}
	\end{figure}

\begin{figure}[!htbp]
		\centering
		\subfigure[Reference solution with $``2+3+0"$ basis.]
		{ \includegraphics[width=0.45\textwidth]{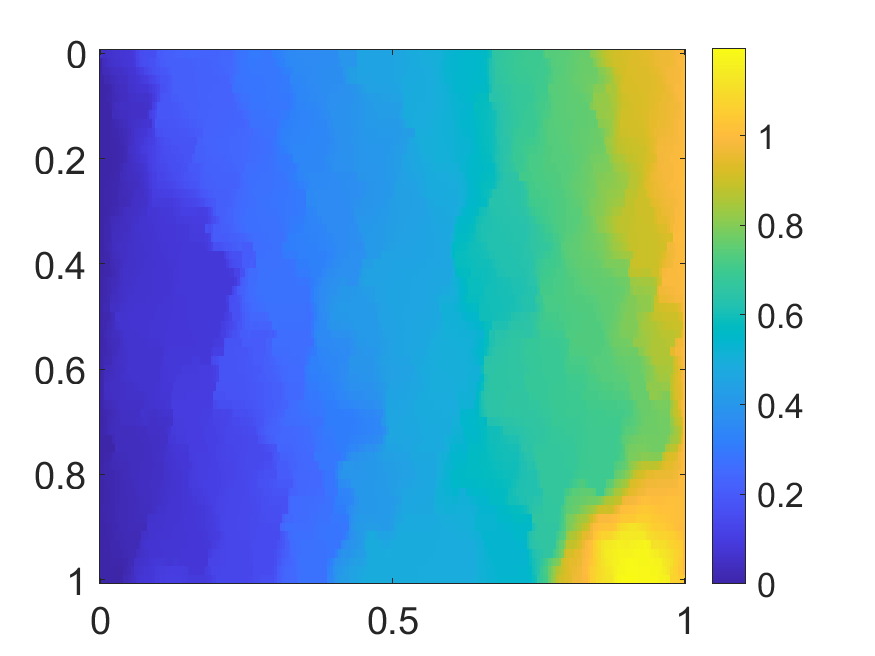}}
		\subfigure[Approximation with $``2+3+0"$ basis.]
		{ \includegraphics[width=0.45\textwidth]{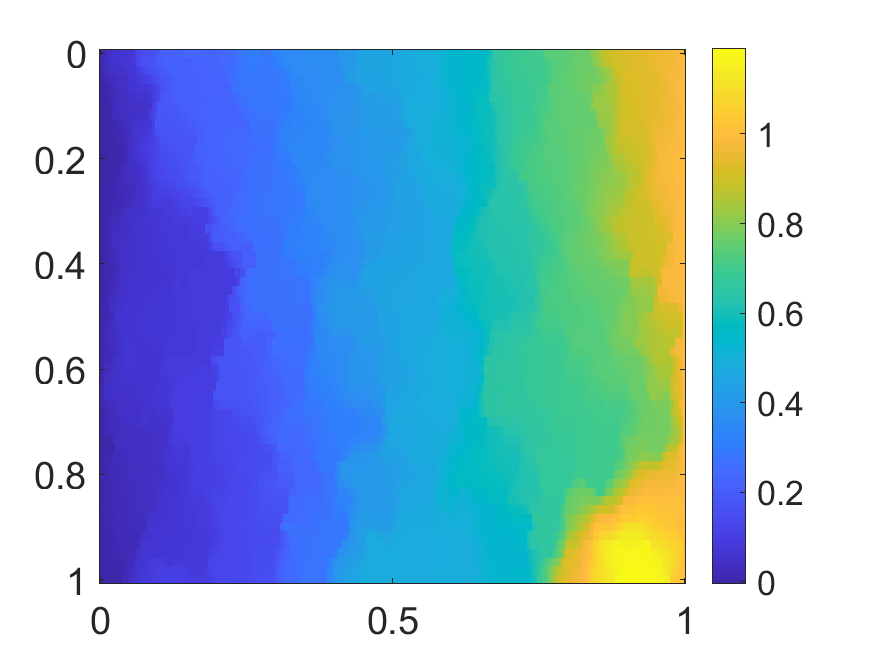}}
		\subfigure[Reference solution with $``2+0+2+1"$ basis.]
		{ \includegraphics[width=0.45\textwidth]{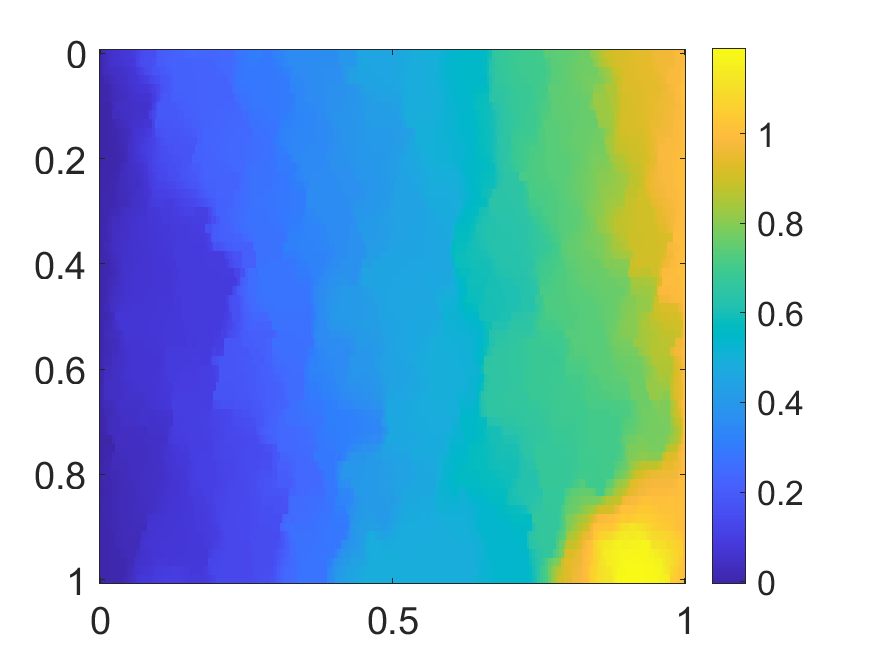}}
		\subfigure[Approximation with $``2+0+2+1"$ basis.]
		{ \includegraphics[width=0.45\textwidth]{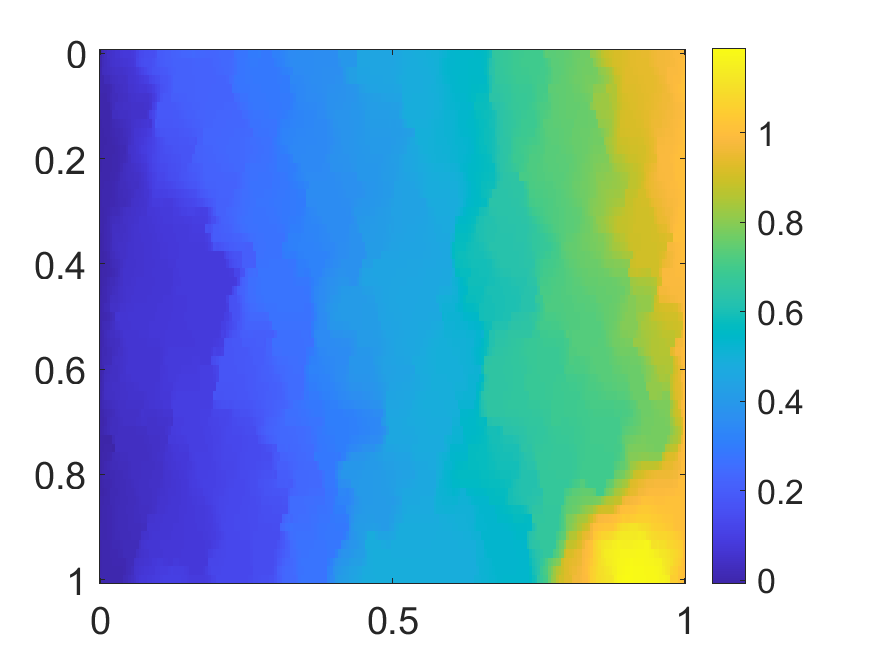}}
		\caption{Comparison of solutions of specific samples at time $t=T$ in step 2 in Algorithm 2. We use SPE model and set $\sigma^2=2$ and $\eta_1=\eta_2=0.5$.}\label{fig:algo2}
	\end{figure}

	We also present the comparison of solutions graphically on some chosen samples, where we use GMsFEM-POD method 1 (Table \ref{tab:algo1}). Here, we use SPE model and set $\sigma^2=2$ and $\eta_1=\eta_2=0.5$.  In Figure \ref{fig:algo1}, we display solutions at final time in three steps. More specifically, for step 1 and 2, we show solutions associated with mean permeability field and some sample fields, where the approximation is GMsFEM solution with $``2+3"$ basis. For step 3, the right column is the POD solution with 10 POD basis. We can have the following conclusions. First, by comparing among three rows in Figure \ref{fig:algo1}, we can see the differences among solutions which are resulted by the uncertainty in permeability field. Secondly, for each row, we could hardly see evident distinctions between the reference and approximation. Specifically, the approximations can relatively well preserve the fine-scale information of the reference solution. In Figure \ref{fig:algo2}, we show comparisons in step 2 using GMsFEM-POD method 1 and 2, respectively. Here, the first row uses algorithm 1 (see Table \ref{tab:algo1}) while the second uses GMsFEM-POD method 2 (see Table \ref{tab:algo2}). We can see both approximations are hardly distinguishable from the reference solution.

	\begin{table}[htbp!] 
		\centering
		\begin{tabular}{c|c|c|c}
			Number of basis & Step1 & Step2 &Step3 \\
			\hline \hline
			2+3  &   1.4\%    &  20\%     &   26\%   \\
			\hline
			5+0  &   30\%   &    28\%   &   34\%   \\
			\hline
			5+5  &   1.7\%      &   9\%   &   23\%  \\
			\hline
			10+0  &  15\%     &   17\%   &26\%      \\
		\end{tabular}
		\caption{Average energy errors (of 100 samples) with different combinations of local basis functions with  Algorithm 1 under $\sigma^2=5$. We set $\eta_1=\eta_2=0.1$ and SPE model.}
		\label{tab:algo1_sigma2_5}
	\end{table}
	\begin{table}[htbp!]
		\centering
		\begin{tabular}{c|c|c|c}
			Number of basis & Step1 & Step2 &Step3 \\
			\hline \hline
			2+3  &   1.4\%    &  8\%     &   9\%   \\
			\hline
			5+0  &   30\%   &  31\%   &   32\%   \\
			\hline
			5+5  &   1.7\% &   5.7\%   &   7\%  \\
			\hline
			10+0  &  15\%     &   16\%   &17\%      \\
		\end{tabular}
		\caption{Average energy errors (of 100 samples) with different combinations of local basis functions with  Algorithm 1 under $\sigma^2=1$. We set $\eta_1=\eta_2=0.1$ and SPE model.}
		\label{tab:algo1_sigma2_1}
	\end{table}
	
	\begin{table}[htbp!]
		\centering
		\begin{tabular}{c|c|c|c}
			Number of basis& Step1 & Step2 & Step3 \\
			\hline \hline
			2+3  &   58\%    & 18.7\%  &19.2\%       \\
			\hline
			2+1+1+1  &  58\% &  18.9\% 13.4\%  8.7\% &9.5\%       \\
			\hline
			5+0  &     30\%  & 30.1\%   &33\%     \\
			\hline
			3+3  & 51\%     &    19.2\%   &21\%   \\
			\hline
			3+1+1+1  & 51\%     & 26\%  18\%  14\%  &15.5\%    \\
			\hline
			3+1+2  & 51\%     & 26\%  15.7\%  &16.5\%    \\
			\hline
			6+0  & 20\%     & 26\%    &27\%\\         
		\end{tabular}
		\caption{Average energy errors (of 100 samples) with different combinations of local basis functions with  Algorithm 2 under $\sigma^2=5$. }
		\label{tab:algo2_sigma2_5}
	\end{table}
	\begin{table}[htbp!]
		\centering
		\begin{tabular}{c|c|c|c}
			Number of basis & Step1 & Step2 & Step3  \\
			\hline \hline
			2+3  &   58\%    & 10.6\%  &11\%       \\
			\hline
			2+1+1+1  &  58\% &  16.2\%  8.3\% 4.9\% & 5.3\%      \\
			\hline
			5+0  &     30\%  & 30.5\%  &32\%      \\
			\hline
			3+3  & 51\%     &    9.4\%  & 10\%   \\
			\hline
			3+1+1+1  & 51\%     & 19.4\%  8.6\%  5.4\%   &6.3\%    \\
			\hline
			3+1+2  & 51\%     & 19.5\%  6.8\%    &7.2\% \\
			\hline
			6+0  & 20\%     & 24\%    &25\% \\   
		\end{tabular}
		\caption{Average energy errors (of 100 samples) with different combinations of local basis functions with  Algorithm 2 under $\sigma^2=1$.}
		\label{tab:algo2_sigma2_1}
	\end{table}
	
	\subsection{Influence of different spatial variability.}
\begin{figure}[htbp!]
		\centering
		\subfigure[Mean of solutions for different $\sigma$]
		{ \includegraphics[width=0.45\textwidth]{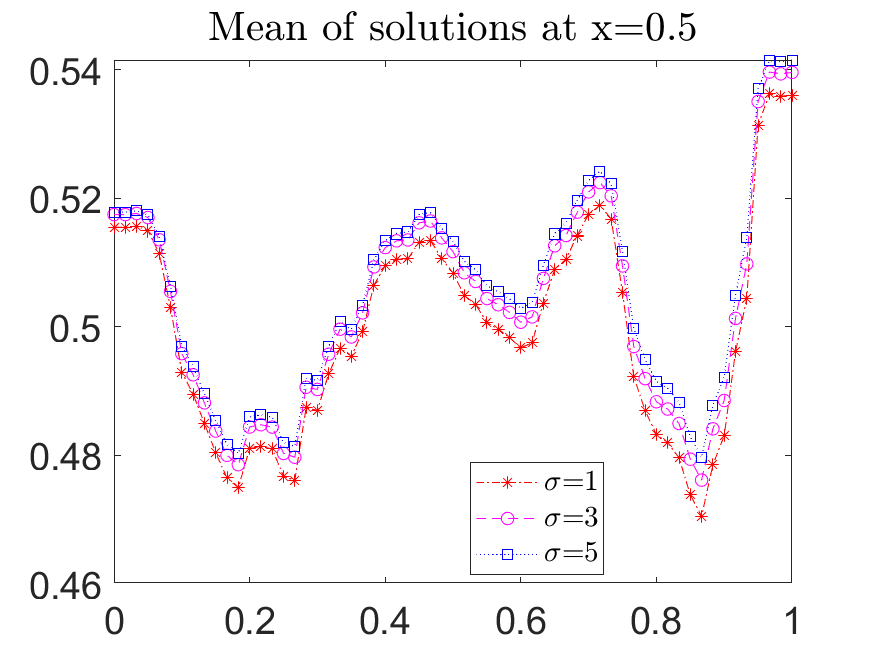}\label{fig:compare_sigma_a}}
		\subfigure[Variance of solutions for different $\sigma$]
		{ \includegraphics[width=0.45\textwidth]{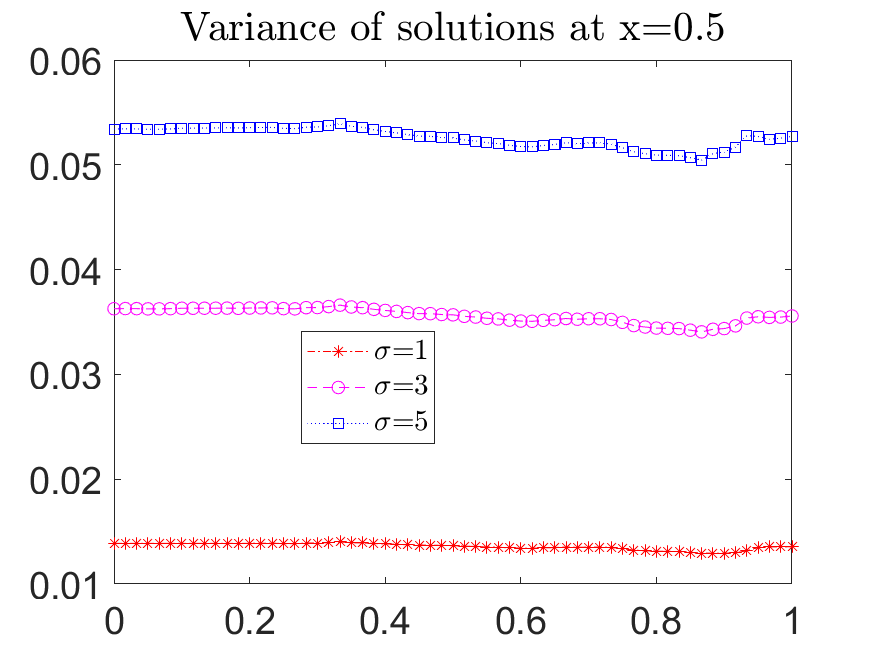}\label{fig:compare_sigma_b}}
		\caption{Comparison of trace of reference solutions  at $x=0.5$ and $t=T$ generated with different $\sigma^2=1,3,5$. Besides, $\eta_1=\eta_2=0.1$ . SPE model (see Figure \ref{fig:model}).}
		\label{fig:compare_sigma}
	\end{figure}

	We also consider the influence of different spatial variability, where we use $\sigma^2=1,3,5$ and $\eta_1=\eta_2=0.1$.  Moreover, we use SPE model. We simply plot the slices of reference solution ( to \eqref{fine sol stoc}) at spatial point $x=0.5$ and $t=T$, i.e. the final time. In Figure \ref{fig:compare_sigma}, for each $\sigma$, we present mean and variance of solutions based on some chosen samples. As one can see from Figure \ref{fig:compare_sigma_a}, there are not apparent differences among three cases. For (b), one can see the variances of solutions are consistent with those of involved random variable $\omega$. Therefore, one can see that the uncertainties in permeability field  indeed influence the solutions in the sense that the intensity of variation in solution  coincides with that of random permeability field. However, it is also shown that the mean solutions are barely affected by the variation of random permeability field.
	\section{Conclusion}
	In this work, we consider a local-global generalized multiscale finite element method for highly heterogeneous stochastic groundwater flow problems. More specifically, we introduce a three-step method (GMsFEM-POD) to obtain a set of basis functions which are independent of the concerned permeability field. In the context of GMsFEM, we use two types of online enrichment methods, adaptive offline enrichment 1 and 2 in Table \ref{online algo1} and \ref{online algo2}) respectively, to construct residual-driven basis functions which are efficient in reducing approximation error. The difference of these two methods exists in whether we keep the residual-driven basis functions obtained in the previous time step. In terms of the permeability field used in performing online enrichment, we further consider two algorithms, GMsFEM-POD method 1 and 2 in Table \ref{tab:algo1} and \ref{tab:algo2}, respectively. In particular, the first one considers mean permeability field to construct residual-driven basis functions while the other uses different sample fields to iteratively perform the online enrichment. We show that the second method can achieve higher accuracy with the same number of basis functions as in the first method. Moreover, we provide analysis for error estimates. The underlying idea of this part is splitting the overall error into four parts and consider each of them individually. In the numerical simulation, we test the proposed methods in three aspects. In particular, we consider the influence of multiscale space dimension with different choices of the parameters in KLE in the first part. In the second part, we compare two offline enrichment algorithms. Lastly, we explore the effect of the spatial variability. Our results can show the efficiency and accuracy of the proposed method, which is also consistent with the analysis part.
	
	
	%

\section*{Acknowledgments}

The research of Eric Chung is partially supported by the Hong Kong RGC General Research Fund (Project numbers 14304719 and 14302018) and the CUHK Faculty of Science Direct Grant 2020-21.

	\bibliographystyle{plain}
	\bibliography{reference}
\end{document}